\documentclass[12pt]{amsart}
\usepackage{amssymb,amsmath,amsthm}  
\usepackage{graphicx}         %
\usepackage{enumerate}
\usepackage{braket}
\usepackage{amsaddr}
\usepackage{amscd}
\usepackage{colordvi}
\usepackage[mathscr]{eucal}
\usepackage{subfigure}
\usepackage{hyperref}
\usepackage{url}
\numberwithin{equation}{section}
\newtheorem{theorem}{Theorem}[section] 

\newtheorem{lemma}[theorem]{Lemma}

\newtheorem{remark}[theorem]{Remark}

\newcommand{\dd}{\displaystyle}
\newcommand{\tu}{\textup}
\newcommand{\bfa}[1]{\boldsymbol{#1}} 			%

\newcommand{\e}{\epsilon}
\newcommand{\ddiv}{\text{div}}     				%
\newcommand{\dx}{ \mathrm{d}x}

\usepackage[numbers]{natbib}
\usepackage{filecontents}
\usepackage[latin9]{inputenc}
\usepackage[letterpaper]{geometry}
\usepackage{pifont}
\usepackage{float}
\usepackage{color}
\usepackage{adjustbox}
\usepackage{diagbox}
\usepackage{multirow}
\usepackage{array}
\usepackage{mathrsfs}
\usepackage{comment}
\usepackage{csquotes}
\graphicspath{{./Figures/}}

\usepackage{amsfonts}

\newcommand{\PH}{\mathcal{T}^H}
\newcommand{\Oh}{\mathcal{T}_h}

\usepackage{color}
\usepackage{adjustbox}
\usepackage{diagbox}
\definecolor{black}{rgb}{0,0,0}

\definecolor{red}{rgb}{1,0,0}

\definecolor{blue}{rgb}{0,0,1}

\usepackage{multirow}

\makeatother

\usepackage[english]{babel}
\author[S\MakeLowercase{hubin} F\MakeLowercase{u}, E\MakeLowercase{ric} C\MakeLowercase{hung}, T\MakeLowercase{ina} M\MakeLowercase{ai}]
{\Large \Large S\MakeLowercase{hubin} F\MakeLowercase{u}, \Large E\MakeLowercase{ric} C\MakeLowercase{hung}, \Large T\MakeLowercase{ina} M\MakeLowercase{ai}$^*$} 
\date{\today}

\title[CEM-GM\MakeLowercase{s}FEM\MakeLowercase{ for nonlinear poroelasticity and elasticity}]
{\textsf{\LARGE C\MakeLowercase{onstraint energy minimizing generalized multiscale 
finite element method for nonlinear poroelasticity and elasticity}}}

\begin{document}

%%%%%%%%%%%%%%%
\begin{abstract}

In this paper, we apply the constraint energy minimizing generalized 
multiscale finite element method (CEM-GMsFEM) to first solving a nonlinear poroelasticity problem.  The arising system consists of a nonlinear pressure equation and a nonlinear stress equation in strain-limiting setting, where strains keep bounded while stresses can grow arbitrarily large.  After time-discretization of the system, to tackle the nonlinearity, we linearize the resulting equations by Picard iteration.  To handle the linearized equations, we employ the CEM-GMsFEM and obtain appropriate offline multiscale basis functions for the pressure and the displacement.  More specifically, first, auxiliary multiscale basis functions are generated by solving local 
spectral problems, via the GMsFEM.  Then, multiscale spaces are constructed in oversampled regions, by solving a constraint energy minimizing (CEM) problem.  After that, this strategy (with the CEM-GMsFEM) is also applied to a static case of the above nonlinear poroelasticity problem, that is, elasticity problem, where the residual based online multiscale basis functions are generated by an adaptive enrichment procedure, to further reduce the error.  Convergence of the two cases is demonstrated by several numerical simulations, which give accurate solutions, with converging coarse-mesh sizes as well as few basis functions (degrees of freedom) and oversampling layers.    

\end{abstract}
%%%%%%%%%%

\maketitle
%%%%%%%%%%%%%%%%


\vskip10pt

%35Q74 PDEs in connection with mechanics of deformable solids
%35J60 Nonlinear elliptic equations
%65N12 Stability and convergence of numerical methods
%65N15 Error bounds
%65N30 Finite elements
%65N99 PDE, BVP 
%74B20 Nonlinear elasticity
%74Q05 Homogenization in equilibrium problems

\vfill

\noindent Shubin Fu $\cdot$ Eric Chung

\noindent Department of Mathematics, The Chinese University of Hong Kong, Shatin, Hong Kong

\noindent E-mail: shubinfu89@gmail.com (Shubin Fu); tschung@math.cuhk.edu.hk (Eric Chung) \\

\noindent Tina Mai\textsuperscript{*} (corresponding author)

\noindent Institute of Research and Development, 
Duy Tan University, Da Nang 550000, Vietnam

\noindent E-mail: maitina@duytan.edu.vn (Tina Mai)

%%%%%

\newpage

\tableofcontents

%%%
\section{Introduction}
For elastic porous media which incompressible viscous fluid flows through, modeling and simulating its deformation are helpful in developing a variety of applications, such as geomechanics or environmental safety.  Given a linear porous medium, Biot \cite{biot} suggested a poroelasticity model, which combines a Darcy flow of the fluid with the behavior of the surrounding linear elastic solid.  In this paper, we investigate a nonlinear poroelasticity model, where the nonlinear stress equation involves quasi-static strain-limiting elasticity (\cite{A-Mai-Walton, B-Mai-Walton}); whereas, the nonlinear pressure equation is a Darcy-type parabolic equation.

To overcome the challenge from the nonlinearity of the poroelasticity, after time-discretization, we use linearization in Picard iteration (with a desired termination criterion) for each time step, until the terminal time.  To tackle the difficulties from multiple scales and high contrast, we apply the constraint energy minimizing generalized 
multiscale finite element method (CEM-GMsFEM, see\cite{cem1,cem2f}) to the linearized equations at the current iteration.  The CEM-GMsFEM here developed from the GMsFEM (\cite{G1}).  

The nonlinear elasticity in the stress equation is motivated by a recent direction of investigating nonlinear responses of materials, thanks to the new developed \textit{implicit constitutive 
theory} (see \cite{Raji03,KRR-ZAMP2007,KRR-MMS2011a,KRR-MMS2011b}).  As Rajagopal remarks, the theory gives a cornerstone to 
developing nonlinear and infinitesimal strain theories for elastic-like (non-dissipative) material behavior.  This setting is different from traditional Cauchy and Green approaches for presenting elasticity which, under the assumption of infinitesimal strains, derive classical linear models.  In addition, it is noteworthy that the implicit constitutive theory yields a stable theoretical base for modeling fluid and solid mechanics diversely, in engineering, chemistry and physics.

Here, in the stress equation, we focus on the \textit{strain-limiting theory} (as a special sub-class of the implicit constitutive theory), where the linearized strain keeps bounded even when the stress becomes extremely large.  
%This contradicts the traditional model (where the stress blows up when the strain blows up and vice versa).
%(In the traditional linear model, the stress blows up when the strain blows up and vice versa.)
Note that it is thus helpful to use the strain-limiting theory to characterize the behavior of fracture, brittle materials near crack tips or notches, or concentrated loads inside the material body (or on its boundary).  Either situation leads to stress intensity despite the small gradient of the displacement (and hence infinitesimal strain).  Within our nonlinear poroelasticity model, the solid part is science-non-fiction and physically valid.  This solid part can undergo infinite stresses and 
does not damage (as the strains are bounded).

Regarding the multiple scales, instead of direct numerical simulations 
on fine grid, model reduction techniques are applied, to lessen the computational burden.  These techniques consist of upscaling and multiscale methods.  On coarse grid, 
upscaling methods mean 
upscaling the material properties based on homogenization, whereas 
multiscale methods need precomputed multiscale 
basis functions.

Within the structure of multiscale methods, in \cite{gnone}, the GMsFEM was used to handle nonlinear problems in poroelasticity.  Then, the idea of CEM-GMsFEM was adopted, for linear poroelasticity in \cite{lporo} (thanks to \cite{dlp}), to create multiscale basis functions (with locally minimal energy) for the pressure and the displacement.  In this paper, for the case of poroelasticity, to deal with the nonlinearity, after the time-discretization, we employ the Picard iteration procedure; and at each iteration, the CEM-GMsFEM is applied as in \cite{lporo}.  The primary component of the CEM-GMsFEM is the construction of local
basis functions for each coarse element (by using the GMsFEM to create the auxiliary multiscale basis
functions) then for each oversampled domain (by employing the CEM to obtain the set multiscale basis
functions).  Convergence analysis within a Picard iteration is shown to support the proposed method.    

As an interesting case of the considered nonlinear poroelasticity, a static strain-limiting nonlinear elasticity model (as in \cite{gne}) is also investigated by similar strategy and via the CEM-GMsFEM.  To take into consideration the influence of source and global information, as in the linear elasticity case (\cite{cemgle}), we use more efficient residual based online basis functions (via adaptive enrichment procedure \cite{chungres}), which are not used in the poroelasticity case (where only offline multiscale basis functions are applied).  The online basis of the CEM-GMsFEM \cite{cem2f} will be computed in an oversampled domain, which is different from the original online approach \cite{chungres}.  We will also provide a proof of global convergence of the Picard iteration procedure in Appendix \ref{conva}. 

Numerical simulations are shown to support the proposed method.  At the end of the Picard iteration process, the CEM-GMsFEM solution is compared with the reference fine-grid solution (at the last time step for the dynamic case). In the static nonlinear elasticity case, we observe that when the sequence of coarse-mesh sizes converges, the sequence of CEM-GMsFEM solutions also accurately converges.  The effects of number of oversampling layers and number of offline multiscale basis functions are as expected.  That is, increasing their numbers (until some certain limits) will increase the CEM-GMsFEM solution accuracy.  The errors further reduce when we adaptively add residual based online basis.  For the nonlinear poroelasticity case, similar conclusions about the CEM-GMsFEM solution (for both the pressure and the displacement) are obtained with respect to the convergence of coarse-grid sizes as well as the oversampling layers.  Regarding the number of offline multiscale basis functions, adding them will improve the displacement accuracy, but will not change the pressure accuracy.            

The next section contains the formulation of our considering 
strain-limiting nonlinear poroelasticity problem.  Section \ref{pre} is for some 
preliminaries about the CEM-GMsFEM, including fine-scale discretization and 
Picard iteration for linearization.  Section \ref{cemforp} is devoted to general idea of the CEM-GMsFEM, for the current nonlinear poroelasticity problem.  Section \ref{cembase} is about computing multiscale spaces, by using the CEM-GMsFEM in our context.  Section \ref{cemfore} discusses an interesting static nonlinear elasticity case of the above nonlinear poroelasticity case.  Numerical results for both cases are provided in Section \ref{nume}.  The last Section \ref{concs} is for conclusions.  In Appendix \ref{conva}, we present a proof of global convergence of the Picard iteration process, by using fixed-point theorem.   
%In Section 7, some numerical results will be shown.

%CEM-GMsFEM \cite{cemgle}

%%%%%%%
\section{Formulation of the nonlinear poroelasticity problem}\label{formulate}
\subsection{Input problem and classical formulation}

Let $\Omega$ be a bounded, Lipschitz, simply connected, open, convex domain of $\mathbb{R}^d \text{ } (d=2,3)$, and $T>0$ be a fixed time.  For the sake of simplicity, the case $d=2$ is considered here.  We refer the readers to our previous paper \cite{gne} for more details about the strain-limiting nonlinear elasticity model.  
We now consider an arising nonlinear poroelasticity system, where the unknowns are displacement $\bfa{u}:  \Omega \times [0,T] $ and pressure $p: \Omega \times [0,T]$ satisfying
\begin{align}
 -\ddiv (\kappa(\bfa{x}, |\bfa{Du}|)\bfa{Du}) + \nabla(\alpha p) &= \bfa{0} \quad \text{in } \Omega \times (0,T]\,,   \label{ue}\\
 \frac{\partial}{\partial t} \left(\alpha \, \ddiv \bfa{u} + \frac{1}{M} p \right) 
 - \ddiv \left( K(\bfa{x}, \bfa{Du},p) \nabla p \right) &= f \quad \text{in }  \Omega \times (0,T]\,, \label{pe}
\end{align}
where the permeability $K(\bfa{x}, \bfa{Du},p)$ can depend on $p$ and $\bfa{Du}$ in non-trivially nonlinear manner (even though our considering materials are isotropic),
%anisotropic (\cite{anper,1diper}), 
 its norm is assumed to be bounded, and
 %Frobenius
\begin{equation}\label{form4k}
  \kappa(\bfa{x}, |\bfa{D}(\bfa{u})|)=\frac{1}{1 - 
  \beta(\bfa{x})|\bfa{D}(\bfa{u})|}\,, 
  %\quad
  %\bfa{a}(\bfa{x},\bfa{D}(\bfa{u})) = 
  %\kappa(\bfa{x}, |\bfa{D}(\bfa{u})|)\bfa{D}%%(\bfa{u})\,,
 \end{equation}
 in which $\bfa{u}(\cdot, t) \in 
   \bfa{W}_0^{1,2}(\Omega)$.  Within this setting, 
   $\kappa(\bfa{x}, |\bfa{D}(\bfa{u})|)\bfa{D}(\bfa{u}) \in 
   \mathbb{L}^1(\Omega)$ and $\bfa{D}(\bfa{u}) \in \mathbb{L}^{\infty}(\Omega)$, as 
   in \cite{Beck2017}.
The boundary and initial conditions are as follows:
\begin{align}
 \bfa{u} & = 0 \quad \text{on }  \partial \Omega \times (0,T]\,, \label{bu}\\
 p &= 0 \quad \text{on } \partial \Omega \times (0,T]\,, \label{bp}\\
 p(\cdot,0) & = p_0 \quad \text{in } \Omega\,, \label{ip}
\end{align}
where $p_0 = 0$ in the numerical simulations (Section \ref{nume}).  To simplify the problem, only homogeneous Dirichlet boundary condition is considered here.  (Other types of boundary conditions can be set simply.)  The heterogeneities are mainly originated from the Cauchy stress tensor $\bfa{T}$, the permeabilities $\kappa$ and $K$, and the Biot-Willis fluid-solid coupling coefficient $\alpha$ (where $\kappa,K, \alpha$ may be highly oscillatory).  We denote by $\nu$ the fluid viscosity and by $M$ the Biot modulus, which are assumed to be constant.  Furthermore, $f$ is a fluid source term (see Theorem \ref{convpt} for its space) representing production or injection processes. 

\begin{remark}\label{npress}
 As an example, another nonlinear poroelasticity problem can be found in \cite{nlss}.  One could consider more general nonlinear form (\cite{gnone}) of $K(\bfa{x}, \bfa{Du},p)$ and use our current Picard linearization technique (as in Section \ref{pre}) to handle the system (\ref{ue})-(\ref{pe}).  Note that our chosen $\kappa(\bfa{x}, |\bfa{Du}|)$ and $K(\bfa{x}, \bfa{Du},p)$ in (\ref{ue})-(\ref{pe}) satisfy the principle of material \textbf{frame-indifference}.
 %as Du symmetric, p const, |Du| const
 Also, for simplicity in our numerical simulations, $K(\bfa{x},\bfa{Du},p)$ can depend only and nonlinearly on $p$ as well as can be a scalar-valued function.  For example, $K(p) = \tu{exp}(p)$ (as in Section 5 in \cite{gnone}).

 %\cite{rajp1}
 %or does not depend on $p$ (as $K$ or $K(x)$ in system (2) of \cite{nslp}) 
\end{remark}

For the stress equation \ref{ue}, in our case of nonlinear elastic stress-strain constitutive relation, the stress tensor $\bfa{T}: \Omega \to \mathbb{R}^{2 \times 2}$ and the traditional \textit{linearized} strain tensor are as follows:
\begin{equation}\label{ted}
 \bfa{T} = \frac{\bfa{D}(\bfa{u})}{1 - \beta(\bfa{x})|\bfa{D}(\bfa{u})|}\,, \qquad 
 \bfa{E} = \bfa{D}(\bfa{u}) = \bfa{Du} =\nabla_s \bfa{u}=\frac{1}{2}(\nabla \bfa{u} + \nabla \bfa{u}^{\text{T}})\,.
\end{equation}
These tensors satisfy our investigating strain-limiting 
model of the following form (\cite{A-Mai-Walton}):
\begin{equation}\label{et}
 \bfa{E} = \frac{\bfa{T}}{1 + \beta(\bfa{x})|\bfa{T}|}\,.
\end{equation}
Equivalently,
\begin{equation}\label{te}
 \bfa{T} = \frac{\bfa{E}}{1 - \beta(\bfa{x})|\bfa{E}|}\,,
\end{equation}
provided that $|\bfa{E}| < \dfrac{1}{\beta(\bfa{x})}$ (which will be explained as follows).

We note that the strain-limiting parameter function $\beta(\bfa{x})$ depends on the position variable $\bfa{x} = (x^1,x^2)$.  From (\ref{et}), it is straightforward that 
\begin{equation}\label{ebound}
 |\bfa{E}| = \frac{|\bfa{T}|}{1 + \beta(\bfa{x})|\bfa{T}|} < \frac{1}{\beta(\bfa{x})}\,,
\end{equation}
which implies that $|\bfa{E}|$ has an upper-bound $\displaystyle \frac{1}{\beta(\bfa{x})}\,.$  Hence, taking large enough $\beta(\bfa{x})$ assures that the limiting-strain owns a small upper-bound, as desired.  
Nevertheless, it is not allowed that $\beta(\bfa{x}) \to \infty$.  %(If $\beta(\bfa{x}) \to \infty$, then $|\bfa{E}| 
%< \displaystyle \frac{1}{\beta(\bfa{x})} 
%\to 0$, a contradiction.)  
Toward the analysis of our problem, $\beta(\bfa{x})$ is assumed to be smooth and possess compact range 
$0 < m_1 \leq \beta(\bfa{x}) \leq m_2\,,$ for some positive constants $m_1, m_2\,.$  Here, we choose $\beta(\bfa{x})$ so that the strong ellipticity condition holds (see \cite{A-Mai-Walton}), that is, $\beta(\bfa{x})$ is sufficiently large, 
to restrain from bifurcations in numerical simulations.

%%%%%
\subsection{Function spaces}
We refer the readers to \cite{C-G-K, gne} for the preliminaries.  Latin indices are in the set $\{1,2\}$.  Functions are denoted by italic capitals (e.g., $f$), vector fields in $\mathbb{R}^2$ and $2 \times 2$ matrix fields over $\Omega$ are denoted by bold letters (e.g., $\bfa{v}$ and $\bfa{T}$).  The space of functions, vector fields in $\mathbb{R}^2$, and $2 \times 2$ matrix fields defined over $\Omega$ are respectively represented by italic capitals (e.g., $L^2(\Omega)$), 
boldface Roman capitals (e.g., $\bfa{V}$), 
and special Roman capitals (e.g., $\mathbb{S}$).  
%The space of symmetric matrices of order 2 is denoted by $\mathbb{S}^2$.  The attached subscript $s$ to a special Roman capital stands for a space of symmetric matrix fields.

Our considering spaces are $\bfa{V}: = \bfa{H}_0^1(\Omega) = \bfa{W}_0^{1,2}(\Omega)$ and $Q: = H^1_0(\Omega)\,.$  The dual norm to $\| \cdot \|_{\bfa{H}_0^1(\Omega)}$ is $\| \cdot \|_{\bfa{H}^{-1}(\Omega)}$.  Here, $| \bfa{v}|$ denotes the Euclidean norm of the 2-component vector-valued function 
$ \bfa{v}$; and $| \nabla \bfa{v}|$ represents the Frobenius norm of the $2 \times 2$ matrix $\nabla \bfa{v}$.  

For every $1 \leq r < 
\infty$, we use $\bfa{L}^r(0,T;\bfa{X})$ to denote the Bochner space with the norm 
\[\|\bfa{w}\|_{\bfa{L}^r(0,T;\bfa{X})} := \left(\int_0^T \|\bfa{w} \|_{\bfa{X}}^r \tu{d}t\right)^{1/r} < + \infty\,, 
\]
\[\|\bfa{w}\|_{\bfa{L}^{\infty}(0,T;\bfa{X})}: = \sup_{0 \leq t \leq T} \|\bfa{w}\|_{\bfa{X}}  < + \infty\,,\]
where $(\bfa{X}, \| \cdot \|_{\bfa{X}})$ is a Banach space.  Also, we define 
\[\bfa{H}^1(0,T; \bfa{X}):= \{ \bfa{v} \in \bfa{L}^2(0,T;\bfa{X}) \, : \, \partial_t \bfa{v} \in \bfa{L}^2(0,T;\bfa{X}) \}\,.\]

Thanks to the notation in \cite{B-M-S}, we will express
\[\bfa{S} \text{ as } \bfa{T} \qquad \text{and} \qquad \bfa{D}(\bfa{u})= \bfa{Du} 
\text{ as } \bfa{E} = \bfa{E}(\bfa{u})\,.\] 

Our current model (\ref{et}) is compatible with the laws of thermodynamics 
\cite{KRR-ARS-PRSA2007,Rajagopal493}, which implies that the class of materials are non-dissipative and elastic. 
   
%%%%%

Thanks to \cite{B-M-S}, 
we derive the following results, which were also stated in 
\cite{BMRS14} (p.\ 19) and proved in our recent GMsFEM paper \cite{gne}.

%%%%% 
\begin{lemma}\label{lem1}
Let 
\begin{equation}\label{calZ}
 \mathcal{Z}:= \left \{ \bfa{\zeta} \in \mathbb{R}^{2 \times 2} \; \biggr | \; 
 0 \leq |\bfa{\zeta}| <\dfrac{1}{m_2}  \right \}\,.
\end{equation}
 For any $\bfa{\xi}\in \mathcal{Z}$ such that $0 \leq |\bfa{\xi}|  
< \dfrac{1}{m_2}$, consider the mapping 
 \[\bfa{\xi} \in \mathcal{Z} \mapsto 
 \bfa{F}(\bfa{\xi}): = \frac{\bfa{\xi}}{1-\beta(\bfa{x}) |\bfa{\xi}|} 
 \in \mathbb{R}^{2 \times 2}\,.\]
 Then, for each $\bfa{\xi}_1, \bfa{\xi}_2 \in \mathcal{Z}$, we have
 \begin{align}
  |\bfa{F}(\bfa{\xi}_1) - \bfa{F}(\bfa{\xi}_2)| 
  &\leq \frac{|\bfa{\xi}_1-\bfa{\xi}_2|} 
  {(1- \beta(\bfa{x})(|\bfa{\xi}_1| + |\bfa{\xi}_2|))^2}\,, \label{cont1} \\
  %if |\bfa{\xi}_2| \leq |\bfa{\xi}_1| then the denominator does not have |\bfa{\xi}_2|.
  (\bfa{F}(\bfa{\xi}_1) - \bfa{F}(\bfa{\xi}_2)) \cdot (\bfa{\xi}_1-\bfa{\xi}_2) 
  &\geq |\bfa{\xi}_1-\bfa{\xi}_2|^2\,. \label{mono1}
 \end{align}
\end{lemma}
%%%

\begin{remark}
 The condition (\ref{mono1}) also means that $\bfa{F}(\bfa{\xi})$ is a monotone operator in 
 $\bfa{\xi}$.
\end{remark}

\begin{remark}\label{caZ}
Without confusion, we will use the condition 
$\bfa{\xi} \in \mathbb{L}^{\infty}(\Omega)$ with the meaning that 
$\bfa{\xi} \in \mathcal{Z}' = \left \{ \bfa{\zeta} \in \mathbb{L}^{\infty}(\Omega) \; \biggr | \; 
 0 \leq |\bfa{\zeta}| <\dfrac{1}{m_2}  \right \}$.
\end{remark}

%Then,
%\begin{equation}\label{dFbound}
%\bfa{0} < D_{\bfa{\xi}}F(\bfa{\xi}) = \frac{1}{(1-\beta(x^1)|\bfa{\xi}|)^2} \leq \left(\frac{1}{1-DM}\right)^2 =C\,.
%\end{equation}
%Recall strong ellipticity formula p.16 paper, give exactly this result.
%Hence, $F(\bfa{\xi})$ is an increasing function, and thus a monotone operator in $\bfa{\xi}$.

%%%%

Let 
\begin{equation}\label{calU}
 \mathcal{U} = \{ \bfa{w} \in \bfa{H}^1(\Omega) \; | 
 \; \bfa{Dw} \in \mathcal{Z}'\}\,,
\end{equation}
with the given $\mathcal{Z}'$ in Remark \ref{caZ}. 
\begin{remark}\label{caU}
Without confusion, we will use the condition 
$\bfa{u}, \bfa{v} \in \bfa{H}_0^1(\Omega) \text{ or } \bfa{H}^1(\Omega)$ 
(context-dependently)
with the meaning that 
$\bfa{u}, \bfa{v} \in \mathcal{U}$.
\end{remark}

%Green formula (integration by parts) 
% \int_{\Omega} T \cdot Dv dx = - \int_{\Omega} (div T) \cdot v dx + \int_{\partial \Omega} (Tn) \cdot v d(\partial \Omega)\,.
% \int u dv = uv - \int v du.
%%%%%
\section{Fine-scale discretization and Picard iteration for linearization}
\label{pre}
We now derive the variational formulation corresponding to the system (\ref{ue})-(\ref{pe}).  First, we multiply Eqs.\ (\ref{ue}) and (\ref{pe}) with test functions from $\bfa{V}$ and $Q$, respectively.  Then, using the Green's formula and the boundary conditions (\ref{bu})-(\ref{ip}), we get the following variational problem: find $\bfa{u}(\cdot,t) \in \bfa{V}$ and $p(\cdot,t) \in Q$ such that
\begin{align}
 a(\bfa{u}, \bfa{v}) - d(\bfa{v},p) &=0\,, \label{vs}\\
 %multiplied with v
 d(\partial_t \bfa{u},q) + c(\partial_t p, q) + b(p,q) &= (f,q) \label{vp}\,,
 %multiplied with q
\end{align}
for all $\bfa{v} \in \bfa{V}$ and $q \in Q$, and the initial pressure is
\begin{equation}\label{vc}
 p(\cdot, 0) = p_0 \in Q\,.
\end{equation}

We define the following nonlinear forms 
\begin{equation}\label{nf}
 a(\bfa{u}, \bfa{v}) = \int_{\Omega} \kappa(\bfa{x}, |\bfa{Du}|) \bfa{Du} \cdot \bfa{Dv} \, \dx\,,
\end{equation}
\begin{equation}\label{nfp}
 b(p,q) = \int_{\Omega} K(\bfa{x}, \bfa{Du},p)  \nabla p \cdot \nabla q \, \dx \,,
\end{equation}
and bilinear and linear forms
\begin{align*}\label{lfdf}
c(p,q) &= \int_{\Omega} \frac{1}{M} p \, q \, \dx \,,\\
 d(\bfa{u},q) &= \int_{\Omega} \alpha (\ddiv \bfa{u}) q \, \dx\,, \qquad (f,q) = \int_{\Omega} f \, q \, \dx\,.
\end{align*}
%Here, $(\cdot , \cdot)$ represents the standard inner product or L2 inner product
Note that (\ref{vs}) can be used to define a relevant initial value $\bfa{u}_0 := \bfa{u}(\cdot , 0) \in \bfa{V}$, provided $p(\cdot,0) = p_0 \in Q$.     

To discretize the variational problem (\ref{vs})-(\ref{vp}), let $\mathcal{T}_h$ (\textit{fine grid}) be a conforming partition for the computational domain $\Omega$, with local grid sizes $h_P: = \tu{diam}(P)\; \forall P \in \mathcal{T}_h$, and $h:= \displaystyle \max_{P \in \mathcal{T}_h} h_P$.  We assume that $h$ is very small so that the fine-scale solution $(\bfa{u}_h, p_h)$ (to be discussed in the following paragraph) is sufficiently near the exact solution.  Next, let $\bfa{V}_h$ and $Q_h$ be the first-order Galerkin (standard) finite element basis spaces with respect to the fine grid $\mathcal{T}_h$, that is,
\[\bfa{V}_h:= \{ \bfa{v} \in \bfa{V}: \bfa{v}|_P \text{ is a polynomial of degree } \leq 1 \; \forall P \in \mathcal{T}_h\}\,,\]
\[Q_h:= \{ q \in Q: q|_P \text{ is a polynomial of degree } \leq 1 \; \forall P \in \mathcal{T}_h\}\,.\]

\textbf{Nonlinear Solve:}  We will first derive the time-discretization of the above system (\ref{vs})-(\ref{vp}), then the nonlinearity will be handled.  

Given an initial pair $(\bfa{u}_0, p_0) \in \bfa{V} \times Q$.  In this section, for simplicity in notation, we will omit the subscript $h$ on the fine grid.  To reach the first goal, we will apply the standard fully implicit (backward Euler) finite-difference scheme (or coupled scheme) 
%backward Euler
for the time-discretization.  It is provided by
\begin{align}
 a(\bfa{u}_{s+1},\bfa{v}) - d(\bfa{v}, p_{s+1}) &= 0\,, \label{dts}\\
 d\left( \frac{\bfa{u}_{s+1} - \bfa{u}_s}{\tau},q\right) + c\left(\frac{p_{s+1} - p_s}{\tau},q\right) + b(p_{s+1},q) = (f_{s+1},q)\,,\label{dtp}
\end{align}
with $\bfa{u}_s = \bfa{u}(\bfa{x}, t_s), \, p_s=p(\bfa{x},t_s), f_{s}=f(t_{s})$, where $t_s = s \tau, \, s= 0,1, \cdots, S, \, S \tau = T$, and $\tau >0$.  Note that $(\bfa{u}_s, p_s)$ represents $(\bfa{u}_{s,h}, p_{s,h})\,.$

After the time-discretization by the fully coupled scheme (\ref{dts})-(\ref{dtp}), we will handle the nonlinearity in space by using a linearization based on Picard iteration.  Indeed, given $(\bfa{u}^n, p^n)$ (which, at the $(s+1)$th time step, represents $(\bfa{u}^n_{s+1,h},p^n_{s+1,h})$) from the previous $n$th Picard iteration step, the nonlinear forms (\ref{nf}) and (\ref{nfp}) at the $(n+1)$th Picard iteration can be respectively linearized as follows:
\[a(\bfa{u}^{n+1},\bfa{v}) \approx a_n(\bfa{u}^{n+1},\bfa{v}):= \int_{\Omega} \kappa(\bfa{x},|\bfa{Du}^n|)\bfa{Du}^{n+1} \cdot \bfa{Dv} \, \dx\,,\]
\[b(p^{n+1},q) \approx b_n(p^{n+1},q):= \int_{\Omega} K(\bfa{x}, \bfa{Du}^n,p^n)  \nabla p^{n+1} \cdot \nabla q \, \dx\,,\]
where
\begin{equation}\label{biform1}
 a_n(\bfa{v}_1,\bfa{v}_2) = \int_{\Omega} \kappa
 (\bfa{x},|\bfa{Du}^n|) 
 \bfa{D} \bfa{v}_1 \cdot \bfa{D} \bfa{v}_2 \, \dx\,, 
\end{equation} 
\begin{equation}\label{biform2}
 b_n(q_1,q_2) = \int_{\Omega} 
 K(\bfa{x}, \bfa{Du}^n,p^n)  
 \nabla q_1 \cdot \nabla q_2 \, \dx \,.
\end{equation}

At the $n$th Picard iteration, the space $\bfa{V}_h$ is equipped with the norm 
\[\|\bfa{v}\|^2_{\bfa{V}_h} = a_n(\bfa{v}, \bfa{v}) \quad \; \forall \; \bfa{v} \in \bfa{V}_h\,,\]
and the space $Q_h$ is equipped with the norm
\[\|q\|^2_{Q_h} = b_n(q, q) \quad \; \forall \; q \in Q_h\,.\]

Provided $(\bfa{u}_s,q_s) \in \bfa{V}_h \times Q_h\,,$ we fix the time-step at $(s+1)$ and take data from the previous Picard iteration $(\bfa{u}^n_{s+1},p^n_{s+1})$ (where we guess a starting point $(\bfa{u}^0_{s+1},p^0_{s+1}) \in \bfa{V}_h \times Q_h$).  For $n=0,1,2, \cdots,$ we wish to find $(\bfa{u}^{n+1}_{s+1}, p^{n+1}_{s+1})$ (that is, $(\bfa{u}^{n+1}_{s+1,h},p^{n+1}_{s+1,h})$) such that
\begin{align}
 a_n(\bfa{u}^{n+1}_{s+1},\bfa{v}) - d(\bfa{v},p^{n+1}_{s+1}) & = 0\,, \label{ls}\\
 d\left( \frac{\bfa{u}^{n+1}_{s+1} - \bfa{u}_s}{\tau},q\right) + c\left( \frac{p^{n+1}_{s+1} - p_s}{\tau},q\right) + b_n(p^{n+1}_{s+1},q) & = (f_{s+1},q)\,.\label{lp}
\end{align}

On the fine grid, the initial value $p_{0,h} \in Q_h$ is set to be the $L^2$ projection of $p_0 \in Q$.  Thus, the initial value $\bfa{u}_{0,h}$ for the displacement is the solution of the equation
\begin{equation}\label{5i}
 a_n(\bfa{u}_{0,h},\bfa{v}) = d(\bfa{v},p_{0,h})\,,
\end{equation}
for all $\bfa{v} \in \bfa{V}_h$. 

We denote by $r$th the Picard iteration where the desired convergence criterion is reached, at the $(s+1)$th time step.  The terminal $(\bfa{u}^{r}_{s+1,h}, p^{r}_{s+1,h})$ now can be set as previous time data, and can be written as $(\bfa{u}_{s+1,h}, p_{s+1,h})$.  

Then, we come back to the algorithm time-stepping (\ref{dts})-(\ref{dtp}) for $s=0,1,\cdots,S$; and within each fixed time, we continue the Picard linearization procedure in (\ref{ls})-(\ref{lp}), until the terminal time $T=S\tau$. 

\begin{remark}
 Theoretically, as in \cite{lporo}, combining Korn's first inequality (\cite{korn}) and the Poincar\'{e} inequality as well as recalling Remark \ref{caZ}, we obtain 
\[ c_T \|\bfa{v}\|_1^2 \leq a_n(\bfa{v}, \bfa{v})=: \|\bfa{v}\|^2_{a_n} \leq C_T \| \bfa{v} \|_1^2\,,
\]
for all $\bfa{v} \in \bfa{V}$, where $c_T$ and $C_T$ are positive constants.  Similarly, there exist two positive constants $c_{\kappa}$ and $C_{\kappa}$ such that 
\[c_{\kappa} \|q\|_1^2 \leq b_n(q,q)=:\|q\|_{b_n}^2 \leq C_{\kappa} \|q\|_1^2\,,\]
for all $q \in Q$.  The existence and uniqueness of solution $(\bfa{u},p)$ for (\ref{ls})-(\ref{lp}) in this linear case can be found in \cite{eulp}. 
\end{remark}

We note that this traditional way will give us a reference fine-scale solution.  The purpose of this paper is to construct a dimension reduction system thanks to (\ref{ls})-(\ref{lp}).  In this spirit, we introduce the reduced finite-dimensional multiscale spaces $\bfa{V}_{\tu{ms}} \subseteq  \bfa{V}, Q_{\tu{ms}} \subseteq Q$, for approximating the solution $(\bfa{u},p)$ on some coarse grid (to lessen the computational cost).

%%%%
\section{CEM-GMsFEM for nonlinear poroelasticity problem}\label{cemforp}
\subsection{Overview}\label{overv}
We will present the construction of auxiliary spaces and multiscale spaces, in the fluid (or pressure) calculation and in the mechanics (or displacement) computation, for the nonlinearly coupled formulation (\ref{vs})-(\ref{vp}).  From the linearized formulation (\ref{ls})-(\ref{lp}), we may view the nonlinearity as constant at each Picard iteration (after time-discretization), to design a suitable CEM-GMsFEM.  In this manner, multiscale spaces are able to be constructed with respect to this nonlinearity.

\textbf{Standard notation.}  Let $\mathcal{T}^H$ be a conforming partition of the domain $\Omega$ such that $\mathcal{T}_h$ is a refinement of $\mathcal{T}^H$.  We call $H:= \displaystyle \max_{K \in \mathcal{T}^H} \tu{diam}(K)$ the coarse-mesh size and $\mathcal{T}^H$ the coarse grid.  Each element of $\mathcal{T}^H$ is called a coarse grid block (element or patch).  
We denote by $N_v$ the total number of interior vertices of 
$\mathcal{T}^H$ and $N$ the total number of coarse blocks (elements).  Let $\{\bfa{x}_i\}^{N_v}_{i=1}$ be the set of vertices (nodes) in $\mathcal{T}^H$ and 
\[w_i= \displaystyle \bigcup_j \Big \{K_j \in \mathcal{T}^H \; \big | \; \bfa{x}_i \in \overline{K_j}\Big \}\]
be the coarse neighborhood of the node $\bfa{x}_i$.  Our main goal is to find a multiscale solution $(\bfa{u}_{\tu{ms}}, p_{\tu{ms}})$ which is a better approximation of the fine-scale solution 
$(\bfa{u}_h, p_h)$ than within GMsFEM (\cite{gnone}).  
%in each Picard iteration and for the convergent ones at the end of Picard iteration procedure.
This is the reason why the CEM-GMsFEM is used to obtain the multiscale solution $(\bfa{u}_{\tu{ms}}, p_{\tu{ms}})$.  

To construct the multiscale spaces, we need two stages.  First, auxiliary spaces are created thanks to the GMsFEM.  Second, using these auxiliary spaces, multiscale spaces are constructed and consist of basis functions whose energy are locally minimized in some subdomains.  After all, these energy-minimized basis functions can be used to obtain a multiscale solution.  

%%%%
\subsection{General idea of the CEM-GMsFEM for nonlinear poroelasticity}\label{gidea}
For details of the GMsFEM and CEM-GMsFEM, we refer the readers to \cite{gne, G1, G2, chungres1, chungres, chung2016adaptive} 
%, chungres1, chungres, chung2016adaptive
and \cite{cem1,cem2f}, respectively.  In this paper, we follow the procedure in Section \ref{pre}, provided $(\bfa{u}_{s,\tu{ms}}, p_{s,\tu{ms}})$ in the multiscale space $\bfa{V}_{\tu{ms}} \times Q_{\tu{ms}} \text{ }(\subset \bfa{V} \times Q)$ (to be discussed later).  At the fixed time $(s+1)$ and current $(n+1)$th Picard iteration, we will use the continuous Galerkin (CG) formulation, with a similar form to the fine-scale problem (\ref{ls})-(\ref{lp}).  More specifically, given the $n$th Picard iteration solution $(\bfa{u}^{n}_{s+1,\tu{ms}}, p^{n}_{s+1,\tu{ms}})$, we wish to find solution $(\bfa{u}^{n+1}_{s+1,\tu{ms}}, p^{n+1}_{s+1,\tu{ms}})$ in $\bfa{V}_{\tu{ms}} \times Q_{\tu{ms}}$ such that
\begin{align}
 a_n(\bfa{u}^{n+1}_{s+1, \tu{ms}},\bfa{v}) - d(\bfa{v},p^{n+1}_{s+1,\tu{ms}}) & = 0\,, \label{lsm}\\
 d\left( \frac{\bfa{u}^{n+1}_{s+1,\tu{ms}} - \bfa{u}_{s,\tu{ms}}}{\tau},q\right) + c\left( \frac{p^{n+1}_{s+1,\tu{ms}} - p_{s,\tu{ms}}}{\tau},q\right) + b_n(p^{n+1}_{s+1,\tu{ms}},q) & = (f_{s+1},q)\,,
 \label{lpm}
\end{align}
for all $(\bfa{v},q) \in \bfa{V}_{\tu{ms}} \times Q_{\tu{ms}}$, with initial condition $p_{0,\tu{ms}} \in Q_{\tu{ms}}$ defined by
\[b(p_{0,h} - p_{0,\tu{ms}}, q) = 0\]
for all $q \in Q_{\tu{ms}}\,.$  The initial value $\bfa{u}_{0,\tu{ms}}$ for the displacement satisfies
\begin{equation}\label{6i}
 a_n(\bfa{u}_{0,\tu{ms}},\bfa{v}) = d(\bfa{v},p_{0,\tu{ms}})\,,
\end{equation}
for all $\bfa{v} \in \bfa{V}_{\tu{ms}}$. 

One notices that the key ingredient of the CEM-GMsFEM is the construction of local basis functions for each coarse element (by applying the GMsFEM to create the auxiliary multiscale basis functions) then for each oversampled domain (by employing the CEM to obtain the multiscale basis functions, which span the multiscale spaces).

%%%%%

\section{Construction of multiscale spaces}\label{cembase}
\begin{figure}
	\centering
	\includegraphics[width=3in]{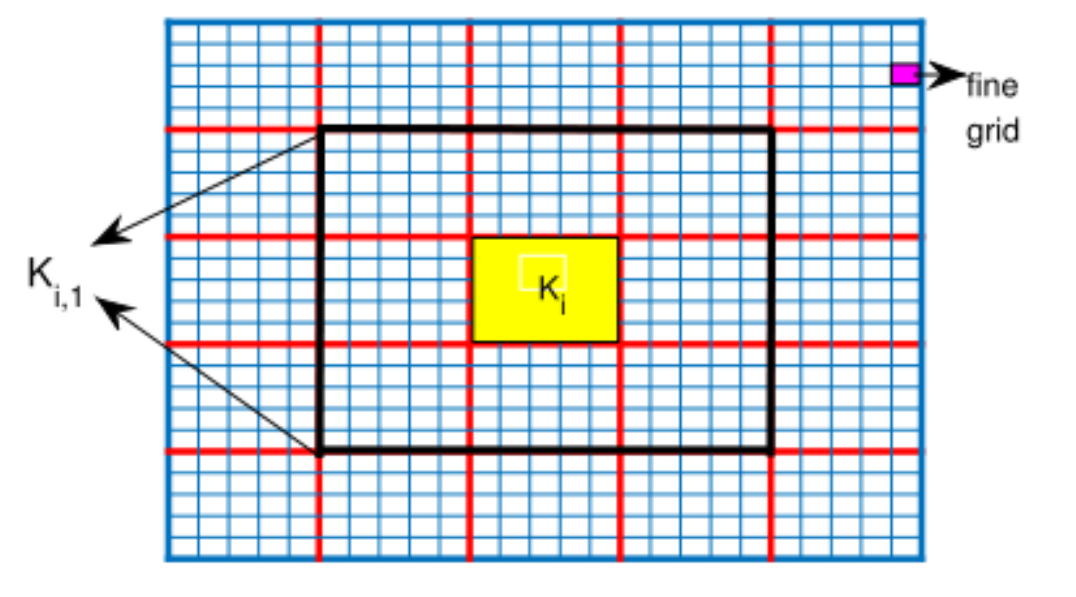}
	\caption{Illustration of the coarse grid $\mathcal{T}^H$, the fine grid $\mathcal{T}_h$, and the oversampled domain $K_{i,1}$.}
	\label{fig:grid}
\end{figure}
This section is devoted to constructing multiscale basis functions, at the  coarse neighborhood $w_i\,,$ with the fixed time $(s+1)$, and the $n$th Picard iteration ($n \geq 0$, given $\bfa{u}^n_{\tu{ms}}:=\bfa{u}^n_{s+1,\tu{ms}} \in \bfa{V}_{\tu{ms}}\,,$ $p^n_{\tu{ms}}:=p^n_{s+1,\tu{ms}} \in Q_{\tu{ms}}$, where the spaces will be explained below).

\subsection{Auxiliary multiscale basis functions}\label{auxb}
We construct auxiliary multiscale basis functions by solving spectral problems on each coarse block $K_i$, making use of the spaces $\bfa{V}(K_i):= \bfa{V}|_{K_i}$, and $Q(K_i):=Q|_{K_i}$.  More specifically, we consider the following local eigenvalue problems: find $(\lambda^i_j, \bfa{v}^i_j) \in \mathbb{R} \times \bfa{V}(K_i)$ such that
\begin{equation}\label{as}
 a_n^i(\bfa{v}^i_j, \bfa{v}) = \lambda_j^i s_n^i(\bfa{v}^i_j, \bfa{v}) \quad \forall \bfa{v} \in \bfa{V}(K_i)\,,
\end{equation}
and find $(\zeta_j^i, q_j^i) \in \mathbb{R} \times Q(K_i)$ such that
\begin{equation}\label{ap}
 b_n^i(q^i_j,q) = \zeta^i_j r_n^i(q^i_j,q) \quad \forall q \in Q(K_i)\,,
\end{equation}
where
\[a_n^i(\bfa{u}, \bfa{v}):= \int_{K_i} (\kappa(\bfa{x},|\bfa{Du}^n_{\tu{ms}}|)\bfa{Du}) \cdot \bfa{Dv} \, \dx\,, \qquad s_n^i(\bfa{u},\bfa{v}):=\int_{K_i} \tilde{\kappa}_s \bfa{u} \cdot \bfa{v} \, \dx\,,\]

\[b_n^i(p,q):= \int_{K_i} K(\bfa{x}, \bfa{Du}^n_{\tu{ms}},p^n_{\tu{ms}}) \nabla p \cdot \nabla q \, \dx\,, \qquad r_n^i(p,q) = \int_{K_i} \tilde{\kappa}_r pq \, \dx\,,\]
in which
\[\tilde{\kappa_s} = \kappa(\bfa{x},|\bfa{Du}^n_{\tu{ms}}|) \sum_{k=1}^{N_v} | \nabla \chi_k^s|^2\,, \qquad \tilde{\kappa_r} = K(\bfa{x}, \bfa{Du}^n_{\tu{ms}},p^n_{\tu{ms}})  \sum_{k=1}^{N_v}  | \nabla \chi_k^r|^2\,.\]
Here, $\chi^s_k$, $\chi^r_k$ are partition of unity functions (\cite{unitybabu}) defined on each neighborhood (that is, for each coarse node) of the coarse mesh (see \cite{gne}, for instance).  More explicitly, for $l=s,r$, the function $\chi_k^l$ satisfies $H|\nabla \chi_k^l| = \mathcal{O}(1), 0 \leq \chi_k^l \leq 1$, and $\displaystyle \sum_{k=1}^{N_v} \chi_k^l =1$.

Assume that the eigenvalues $\{\lambda_j^i\}$ as well as $\{\zeta_j^i\}$ are ordered ascendingly, and the eigenfunctions satisfy the normalization condition $s_n^i(\bfa{v}^i_j, \bfa{v}^i_j) = 1$ as well as $r_n^i(q^i_j,q^i_j) = 1$.  Next, we pick $J_i^v \in \mathbb{N}^+$ and define the local auxiliary space $\bfa{V}_{\tu{aux}}(K_i):=\tu{span} \{\bfa{v}_j^i: 1 \leq j \leq J_i^v \}$.  In the same way, we choose $J_i^q \in \mathbb{N}^+$ and define $Q_{\tu{aux}}(K_i): = \tu{span} \{q_j^i: 1 \leq j \leq J_i^q\}$.  Thanks to these local spaces, we define the global auxiliary spaces $\bfa{V}_{\tu{aux}}$ and $Q_{\tu{aux}}$ by 
\begin{align*}
\bfa{V}_{\tu{aux}}:=  \bigoplus_{i=1}^{N} \bfa{V}_{\tu{aux}}(K_i) \subseteq \bfa{V} \, \text{  and  } Q_{\tu{aux}}:=  \bigoplus_{i=1}^{N} Q_{\tu{aux}}(K_i) \subseteq Q\,.
\end{align*}

The inner products of the global auxiliary multiscale spaces are defined by
\begin{align*}
s_n(\bfa{u}, \bfa{v}):= \sum_{i=1}^{N} s_n^i(\bfa{u},\bfa{v})\,, \qquad \| \bfa{v} \|_{s_n}:= \sqrt{s_n(\bfa{v}, \bfa{v})}\, \quad \forall \bfa{u}, \bfa{v} \in \bfa{V}_{\tu{aux}}\,, \\
r_n(p,q):= \sum_{i=1}^{N} r_n^i(p,q)\,, \qquad \|q\|_{r_n}:= \sqrt{r_n(q,q)}\, \quad \forall p,q \in Q_{\tu{aux}}\,.
\end{align*}

Moreover, defining projection operators $\pi_n^v: \bfa{V} \to \bfa{V}_{\tu{aux}}$ and $\pi_n^q: Q \to Q_{\tu{aux}}$ such that for all $\bfa{v} \in \bfa{V}, q \in Q$, we have 
\begin{align*}
 \pi_n^v(\bfa{v}):= \sum_{i=1}^{N} \sum_{j=1}^{J_i^v} s_n^i(\bfa{v},\bfa{v}^i_j)\bfa{v}^i_j\,, \quad \pi_n^q(q):= \sum_{i=1}^{N} \sum_{j=1}^{J_i^q} r_n^i(q,q^i_j)q^i_j\,.
\end{align*}

%%%%%
\subsection{Multiscale spaces}\label{mss}
Now, we construct the multiscale spaces toward the practical simulations.  For each coarse block $K_i$, we define the oversampled subdomain $K_{i,m} \subset \Omega$ by expanding $K_i$ by $m$ layers, that is,
\[K_{i,0}:=K_i, \quad K_{i,m} :=\bigcup \Big \{K \in \mathcal{T}^H: K \cap K_{i,m-1} \neq
%overline{....}
\emptyset \Big \}, \quad m = 1,2, \cdots\,.\]
%See Fig.\ \ref{fig:grid} for illustration of the coarse grid $\mathcal{T}^H$, the fine grid $\mathcal{T}_h$, and the oversampled region $K_{i,1}\,.$  
We define    
\[\bfa{V}(K_{i,m}):= [\bfa{H}_0^1(K_{i,m})]^d\,, \quad Q(K_{i,m}):=H_0^1(K_{i,m})\,.\] 
Here, $d=2\,.$  

\begin{figure}[H]
	\centering
	\subfigure[first component of $\bfa{\psi}_{1}^{i,m}$ in $\bfa{V}_{\text{ms}}$]{
		\includegraphics[width=2.in]{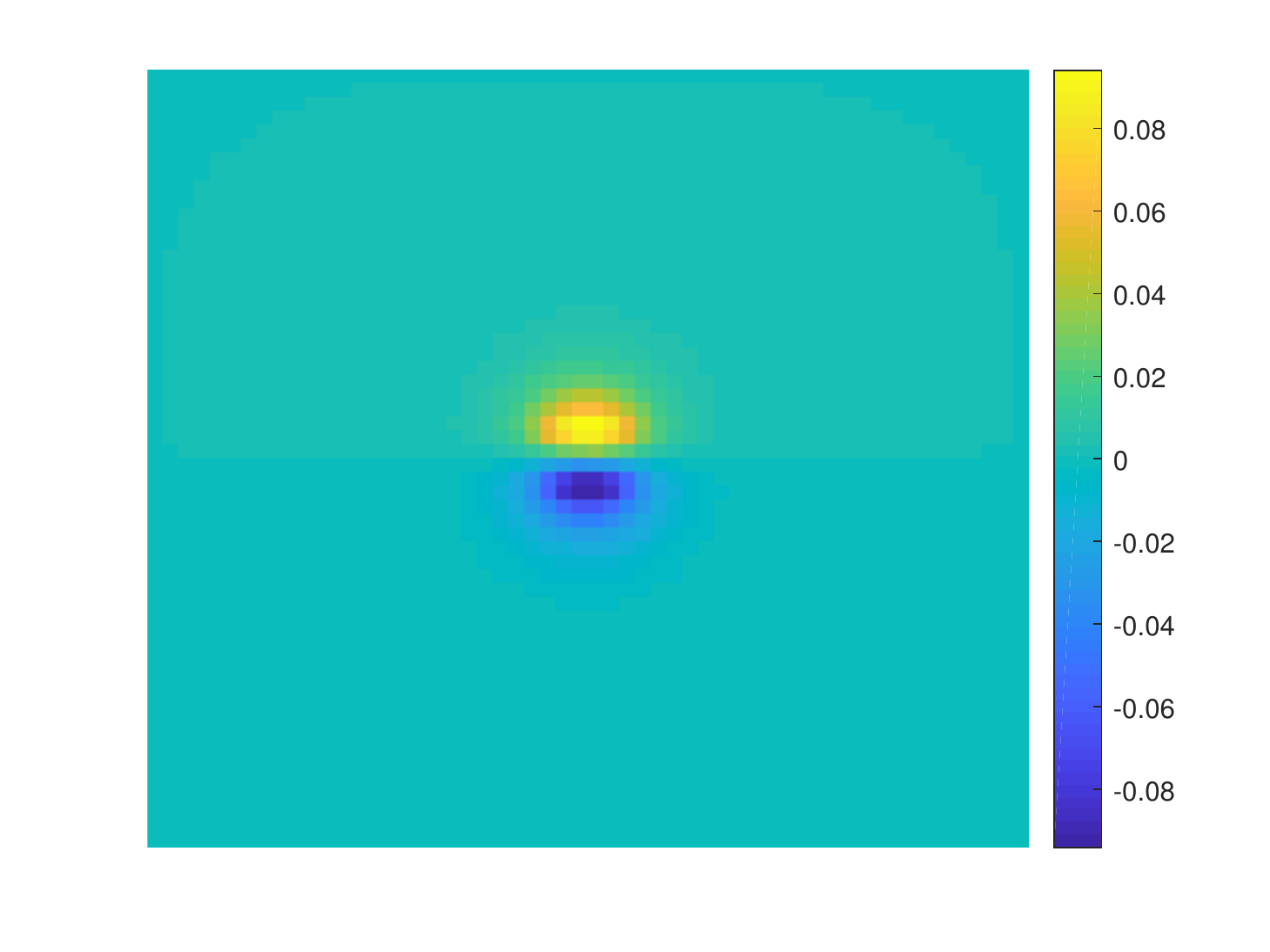}}
	\subfigure[second component of $\bfa{\psi}_{1}^{i,m}$ in $\bfa{V}_{\text{ms}}$]{
		\includegraphics[width=2.in]{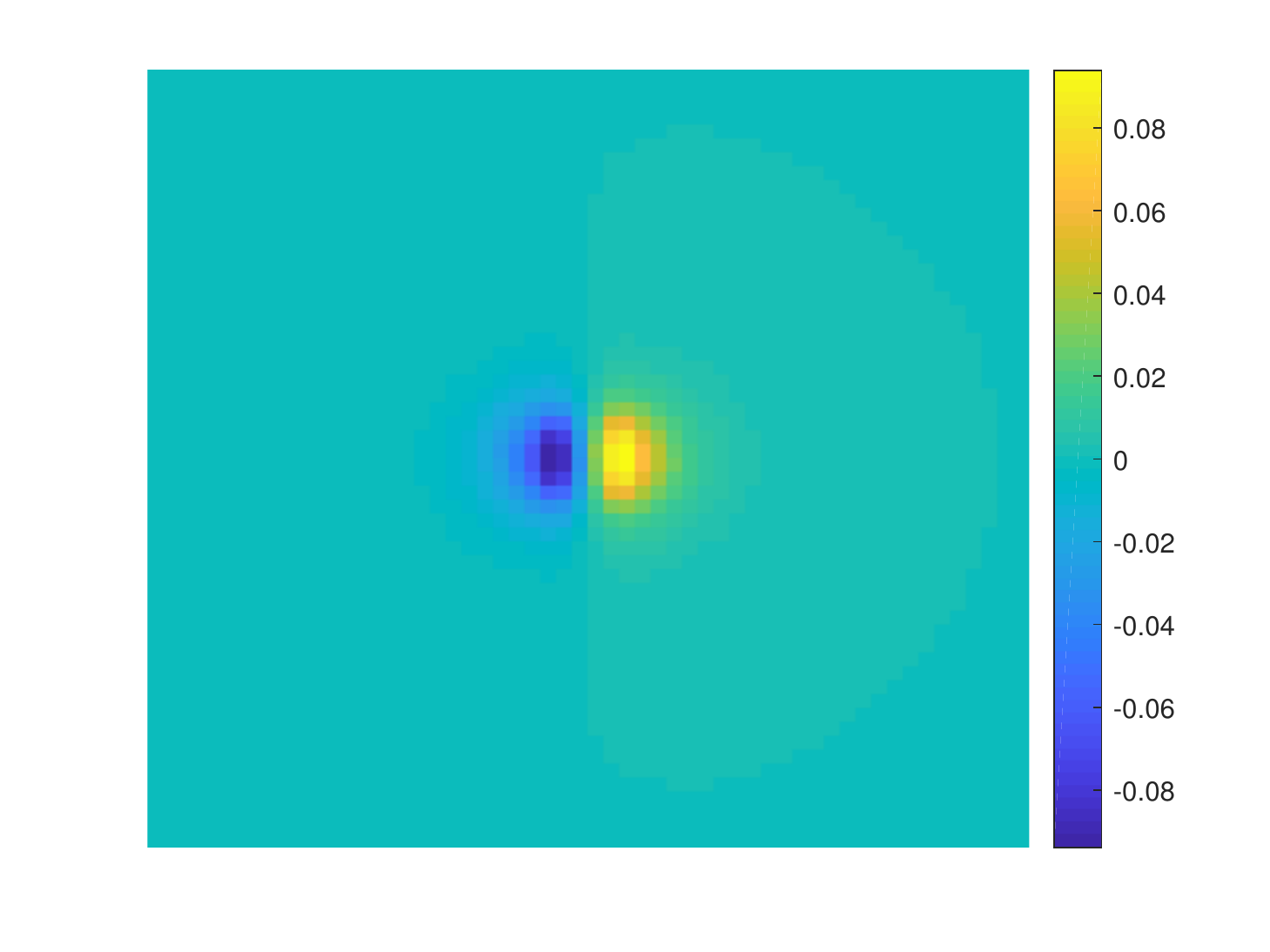}}
	\subfigure[$\phi_{1}^{i,m}$ in $Q_{\text{ms}}$]{
		\includegraphics[width=2.in]{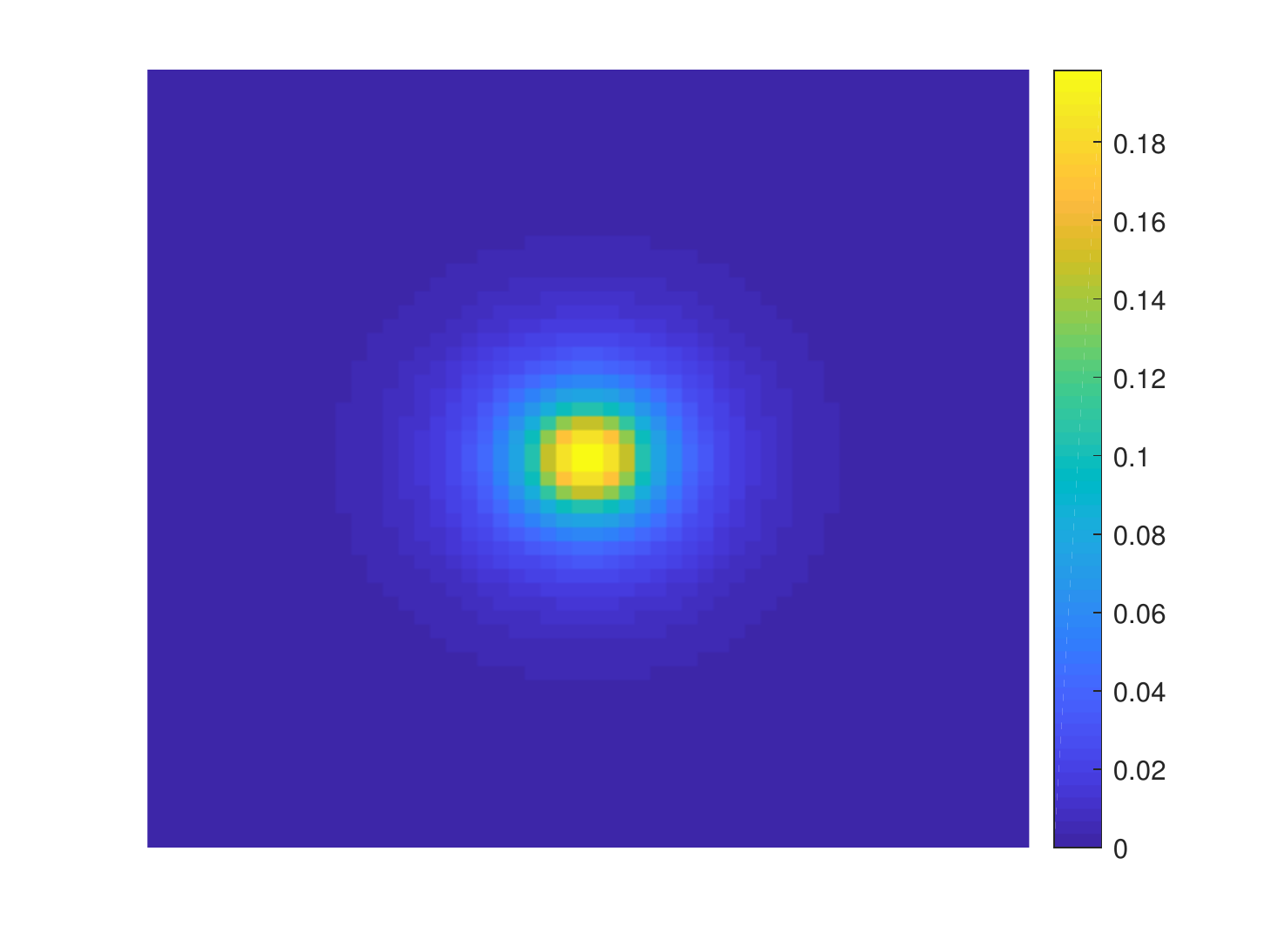}}
	\caption{First multiscale basis functions of $\bfa V_{\text{ms}}$ and $Q_{\text{ms}}$, respectively, for $m=5$ and $H =\sqrt{2}/40$.}
	\label{fig:basis} 
\end{figure}

 After that, for every pair of auxiliary functions $\bfa{v}^i_j \in \bfa{V}_{\tu{aux}}$ and $q_j^i \in Q_{\tu{aux}}$, we solve the following minimization problems: find multiscale basis function $\bfa{\psi}^{i,m}_j \in \bfa{V}(K_{i,m})$ such that 
\begin{equation}\label{cs}
\bfa{\psi}^{i,m}_j = \tu{argmin} \{a_n(\bfa{\psi},\bfa{\psi}) + s_n(\pi^v_n(\bfa{\psi}) -\bfa{v}^i_j, \pi_n^v(\bfa{\psi}) -\bfa{v}^i_j): \bfa{\psi} \in \bfa{V}(K_{i,m}) \} 
 \end{equation}
and find $\phi^{i,m}_j \in Q(K_{i,m})$ such that 
\begin{equation}\label{cp}
\phi^{i,m}_j = \tu{argmin} \{b_n(\phi,\phi) + r_n(\pi_n^q(\phi) - q_j^i, \pi_n^q(\phi) - q_j^i): \phi \in Q(K_{i,m}) \}\,.
 \end{equation}
We note here that the problem (\ref{cs}) is equivalent to the local variational problem
\[a_n\left(\bfa{\psi}^{i,m}_j, \bfa{v}\right) + s_n\left(\pi^v_n(\bfa{\psi}^{i,m}_j), \pi^v_n(\bfa{v})\right) = s_n\left(\bfa{v}^i_j, \pi_n^v(\bfa{v})\right), \quad \forall \bfa{v} \in \bfa{V}\left(K_{i,m}\right)\,,\]
while the problem (\ref{cp}) is equivalent to 
\[b_n\left(\phi_j^{i,m},q\right) + r_n\left(\pi_n^q(\phi^{i,m}_j),\pi_n^q(\bfa{v})\right) = r_n\left(q_j^i, \pi_n^q(q)\right), \quad \forall q \in Q\left(K_{i,m}\right)\,.\]
Last, for fixed parameters $m, J_i^v,J_i^q$, the multiscale spaces $\bfa{V}_{\tu{ms}}$ and $Q_{\tu{ms}}$ are defined through
\[\bfa{V}_{\tu{ms}} := \tu{span} \{ \bfa{\psi}^{i,m}_j: 1 \leq j \leq J_i^v, 1 \leq i \leq N \}\,,\]
and 
\[Q_{\tu{ms}}: = \tu{span} \{ \phi^{i,m}_j: 1 \leq j \leq J_i^q, 1 \leq i \leq N\}\,.\]
See Fig.\ \ref{fig:basis} for illustration of multiscale basis functions.
Similarly, we can interpret the multiscale basis functions $\bfa{\psi}^{i,m}_j \in \bfa V_{\text{ms}}$ and $\phi_j^{i,m} \in Q_{\text{ms}}$ as approximations to global multiscale basis functions $\bfa{\psi}_j^i \in \bfa{V}$ and $\phi_j^i \in Q$ by
\begin{equation}\label{csg}
\bfa{\psi}^{i}_j = \tu{argmin} \{a_n(\bfa{\psi},\bfa{\psi}) +  s_n(\pi^v_n(\bfa{\psi}) -\bfa{v}^i_j, \pi_n^v(\bfa{\psi}) -\bfa{v}^i_j): \bfa{\psi} \in \bfa{V} \} 
 \end{equation}
and  
\begin{equation}\label{cpg}
\phi^{i}_j = \tu{argmin} \{b_n(\phi,\phi) + r_n(\pi_n^q(\phi) - q_j^i, \pi_n^q(\phi) - q_j^i): \phi \in Q \}\,.
 \end{equation}

 These basis functions are globally supported in the domain $\Omega$, but exponentially decay (as shown in \cite{cem1}) outside some local (oversampled) subdomain.  This feature takes a crucial part in the convergence analysis of the CEM-GMsFEM and proves the use of local multiscale basis functions in $\bfa{V}_{\tu{ms}}$ and $Q_{\tu{ms}}$ (\cite{lporo}).  

%%%%
\subsection{Multiscale method}\label{msm}
In the previous Subsections \ref{auxb} and \ref{mss}, the spaces $\bfa{V}$ and $Q$ are continuous.  Toward computations, we need some finite dimensional analogues of the multiscale spaces $\bfa{V}_{\tu{ms}}$ and $Q_{\tu{ms}}$.  Thus, in our numerical simulations, we solve the considered problem using the fine mesh defined in $K_{i,m}$, via an appropriate finite element method (\cite{cem2f}).  
%after (5), Subsection 3.2

Given $\bfa{u}_{s,\tu{ms}}$ and $p_{s,\tu{ms}}\,,$ fixing the time-step at $(s+1)$, we thus have the following fully discrete scheme for the Picard iteration procedure: choose a starting guess of $\bfa{u}^{\tu{old}}_{s+1,\tu{ms}}$, and compute the multiscale space $\bfa{V}^{\tu{old}}_{s+1, \tu{ms}}$; we then wish to find $(\bfa{u}^{\tu{new}}_{s+1,\tu{ms}}, p^{\tu{new}}_{s+1,\tu{ms}}) \in \bfa{V}^{\tu{old}}_{s+1,\tu{ms}} \times Q^{\tu{old}}_{s+1,\tu{ms}}$ such that
\begin{align}
 a_{\tu{old}}(\bfa{u}^{\tu{new}}_{s+1,\tu{ms}},\bfa{v}) - d(\bfa{v},p^{\tu{new}}_{s+1, \tu{ms}}) & = 0\,, \label{ds}\\
 d\left( \frac{\bfa{u}^{\tu{new}}_{s+1, \tu{ms}} - \bfa{u}_{s, \tu{ms}}}{\tau},q\right) + c\left( \frac{p^{\tu{new}}_{s+1, \tu{ms}} - p_{s,\tu{ms}}}{\tau},q\right) + b_{\tu{old}}(p^{\tu{new}}_{s+1, \tu{ms}},q) & = (f_{s+1},q)\,.\label{dp}
\end{align}
for all $(\bfa{v},q) \in \bfa{V}^{\tu{old}}_{s+1,\tu{ms}} \times Q^{\tu{old}}_{s+1,\tu{ms}}$ with initial condition $p_{0,\tu{ms}} \in Q^{\tu{old}}_{\tu{ms}}$ defined by
\[b_{\tu{old}}(p_{0,h}  - p_{0,\tu{ms}}, q) = 0\,,\]
for all $q \in Q^{\tu{old}}_{s+1,\tu{ms}}$.  The initial value $\bfa{u}_{0,\tu{ms}}$ for the displacement is the solution of the equation from (\ref{vs}):
\begin{equation}\label{7i}
 a_{\tu{old}}(\bfa{u}_{0,\tu{ms}},\bfa{v}) = d(\bfa{v},p_{0,\tu{ms}})\,,
\end{equation}
for all $\bfa{v} \in \bfa{V}^{\tu{old}}_{s+1,\tu{ms}}$.
Again, we use $r$th to denote the Picard iteration where the desired convergence criterion is reached.  The terminal $(\bfa{u}^{r}_{s+1,\tu{ms}}, p^{r}_{s+1,\tu{ms}})$ can be set as previous time data, and can be written as $(\bfa{u}_{s+1,\tu{ms}}, p_{s+1,\tu{ms}})$.  

We then return to the algorithm time-discretization (\ref{dts})-(\ref{dtp}) for $s=0, 1, \cdots,S$; and continue the iterative Picard linearization (\ref{ds})-(\ref{dp}) until the terminal time $T= S\tau$.

\bigskip

Finally, we derive some convergence result of the CEM-GMsFEM for this dynamic case.  At the time step $s \text{ }( 1\leq s \leq S)$ defined in (\ref{dts})-(\ref{dtp}), within the $n$th Picard iteration ($n \geq 1$), the following result and its proof are obtained directly from \cite{lporo}.
\begin{theorem}\label{convpt}
 Assume sufficiently large parameters $m,J_i^v,J_i^q$, a source function $f \in L^{\infty}(0,T;L^2(\Omega)) \cap H^1(0,T;H^{-1}(\Omega))$ as well as initial data $p_{0,h} \in Q_h$ and $\bfa{u}_{0,h} \in \bfa{V}_h$ defined in (\ref{5i}).  Then, the error between the multiscale solution $(\bfa{u}^n_{s,\tu{ms}}, p^n_{s,\tu{ms}}) \in \bfa{V}_{\tu{ms}} \times Q_{\tu{ms}}$ of (\ref{ds})-(\ref{dp}) and the fine-scale solution $(\bfa{u}^n_{s,h}, p^n_{s,h}) \in \bfa{V}_h \times Q_h$ of (\ref{ls})-(\ref{lp}) satisfies
 \[\|\bfa{u}^n_{s,h} - \bfa{u}^n_{s,\tu{ms}}\|_1 + \| p^n_{s,h} - p^n_{s, \tu{ms}}\|_1 \lesssim H\mathcal{W}_s + t_s^{-1/2} H \|p_{0,h}\|_1\]
 for $s=1,2,\cdots,S$.  Here, $\mathcal{W}_s$ only depends on the data and is defined through
 \[\mathcal{W}_s:= \|p_{0,h} \|_1 + \| f\|_{L^2(0,t_s;L^2(\Omega))}+ \|f\|_{L^{\infty}(0,t_s;L^2(\Omega))} + \|\partial_t f\|_{L^2(0,t_s;H^{-1}(\Omega))}\,.\]
\end{theorem}

\bigskip

%%%%

\section{CEM-GMsFEM for nonlinear elasticity problem}\label{cemfore}

We now consider a static case of the above nonlinear poroelasticity case (from Section \ref{cemforp}), namely nonlinear elasticity problem.
%%%%%from elasticity file%%%
\subsection{Formulation of the problem}\label{formulatee}
\subsubsection{Input problem and classical formulation}\label{input}

We refer the readers to our previous paper \cite{gne} for more details.  Here, we briefly introduce the formulation.  Let our computational domain be $\Omega \in \mathbb{R}^2$ (as in Section \ref{formulate}), which is a strain-limiting nonlinear elastic composite material.  

The material is assumed to be at a static state (\cite{B-M-S}) after 
the action of body forces $\bfa{f}: \Omega \to \mathbb{R}^2$ and traction forces 
$\bfa{G}: \partial \Omega_T 
\to \mathbb{R}^2$.  We denote the boundary of the set $\Omega$ by $\partial \Omega$, 
which is Lipschitz continuous,  
having two parts $\partial \Omega_T$ and $\partial \Omega_D$, with the given displacement 
$\bfa{u}: \Omega \to \mathbb{R}^2$ on $\partial \Omega_D$.  
We are investigating the strain-limiting 
model (\cite{gne}) in either physical form (\ref{et}) or its equivalent mathematical form (\ref{ted}), that is
%(\ref{ted}), (\ref{et}), (\ref{te}), or 
\begin{equation}\label{tea}
 \bfa{T} = \frac{\bfa{D}(\bfa{u})}{1 - \beta(\bfa{x})|\bfa{D}(\bfa{u})|}\,.
\end{equation}

%%%%%
\subsubsection{Function spaces}
We refer the readers to \cite{gne, C-G-K} for the preliminaries, and to Section \ref{formulate} for function spaces.  Let 
%and let $\bfa{u} \in \bfa{V}$ (thus $\bfa{E} \in \mathbb{L}^2(\Omega)$),
\begin{equation}\label{H1*}
 \bfa{f} \in \bfa{H}^{1}_*(\Omega)= 
\left \{\bfa{g} \in \bfa{H}^1(\Omega) \biggr | \int_{\Omega} \bfa{g} \, \dx = \bfa{0}\right \} 
\subset \bfa{L}^2(\Omega) \subsetneq \bfa{H}^{-1}(\Omega)
\end{equation}  
%H_0^1 \subset H^1 \subset L^2 \subset H_0^1 (the last is strict)
be bounded in $\bfa{L}^2(\Omega)$.  The problem we are considering is as follows:  
find $\bfa{u} \in \bfa{H}^1(\Omega)$ and $\bfa{T} \in \mathbb{L}^1(\Omega)$ such that
\begin{align}\label{form1a}
\begin{split}
 -\textup{\ddiv} (\bfa{T}) &= \bfa{f}  \quad \text{in } \Omega \,,\\
  \bfa{Du} &= \frac{\bfa{T}}{1 + \beta(\bfa{x})|\bfa{T}|} \quad \text{in } 
 \Omega \,,\\
 \bfa{u} &= \bfa{0} \quad \text{on } \partial \Omega_D\,,\\
 %only homogeneous BC is considered here
 \bfa{Tn} &= \bfa{G} \quad \text{on } \partial \Omega_T\,,
 %later, we assume $\partial \Omega_T = \emptyset
 \end{split}
 \end{align}
where $\bfa{n}$ denotes the outer unit normal vector to the boundary of $\Omega$.  

Our current model (\ref{et}) is compatible with the laws of thermodynamics 
\cite{KRR-ARS-PRSA2007,Rajagopal493}, which implies that the studying class of materials is non-dissipative and elastic. 

Assuming that $\partial \Omega_T = \emptyset$, we consider an interesting static case of the problem (\ref{ue})-(\ref{pe}), that is, we investigate the following displacement problem of (\ref{form1a}):  
   Find $\bfa{u} \in \bfa{H}_0^1(\Omega)$ such that
  \begin{align}
 -\textup{\ddiv} \left(\frac{\bfa{D}(\bfa{u})}{1 - \beta(\bfa{x})|\bfa{D}(\bfa{u})|}\right) 
 &= \bfa{f} \quad \text{in } 
 \Omega \,,\label{form3} \\ 
 \bfa{u}&= \bfa{0}
 \quad \text{on } \partial \Omega\,.\label{D}
 %\bfa{S}_{\epsilon} = G, \quad (S_{\epsilon})_{(1)} &= (S_{\epsilon})_{(2)} 
 %\quad \text{on } \partial \Omega_S\,,\label{T}
 \end{align}
 Given $\kappa$ from (\ref{form4k}), we denote 
 \begin{equation}\label{form4n}
  %\kappa(\bfa{x}, |\bfa{D}(\bfa{u})|)=\frac{1}{1 - 
  %\beta(\bfa{x})|\bfa{D}(\bfa{u})|}\,, \quad
  \bfa{a}(\bfa{x},\bfa{D}(\bfa{u})) = 
  \kappa(\bfa{x}, |\bfa{D}(\bfa{u})|)\bfa{D}(\bfa{u})\,,
 \end{equation}
 in which $\bfa{u}(\bfa{x}) \in 
   \bfa{W}_0^{1,2}(\Omega)$.  Within this setting, 
   $\bfa{a(\bfa{x},\bfa{\xi})} \in 
   \mathbb{L}^1(\Omega)$, $\bfa{\xi} \in \mathbb{L}^{\infty}(\Omega)$, as 
   in \cite{gne, Beck2017}.  
   
   %Note that the problem (\ref{form3})-(\ref{D}) arises from a strain-limiting model (the details are in \cite{gne}).
   % which will be explained in later paragraph.
   
%%%%%
\subsubsection{Existence and uniqueness}
   For $\bfa{u} \in \bfa{V} =\bfa{H}_0^1(\Omega)$, we multiply Eq.\ (\ref{form3}) 
by $\bfa{v} \in \bfa{V}$ and integrate the resulting equation 
with respect to $\bfa{x}$ over $\Omega$.  
Integrating the first term by parts and using the condition 
$\bfa{v} = \bfa{0}$ on $\partial \Omega$, we obtain
%\cite Short standard in PDEs existence
\begin{equation}\label{w8.2}
 \int_{\Omega} \bfa{a}  (\bfa{x},\bfa{D}\bfa{u} ) \cdot \bfa{D}\bfa{v} \, \dx = 
 \int_{\Omega} \bfa{f} \cdot \bfa{v} \, \dx 
 \,, \quad \forall \bfa{v} \in \bfa{V} \,.
\end{equation}
By the weak (often called generalized) formulation of the boundary value problem 
(\ref{form3})-(\ref{D}), we interpret the problem as follows:
\begin{equation}\label{w8.3}
\text{ Find } (\bfa{u}, \bfa{Du})  \in \bfa{V} \times 
\mathbb{L}^{\infty}(\Omega), 
\text{that is, find } 
\bfa{u} \in \bfa{V}
\text{ such that } (\ref{w8.2}) \text{ holds for each } \bfa{v} \in \bfa{V}\,.
\end{equation} 

We refer the readers to our previous paper \cite{gne} and references therein for the existing results about the existence and uniqueness of the weak solution $\bfa{u} \in \bfa{H}_0^1(\Omega)$ to (\ref{w8.3}), or $(\bfa{u}, \bfa{T}) \in  \bfa{H}_0^1(\Omega) \times \mathbb{L}^1(\Omega) \text{ (or } \bfa{H}_0^1(\Omega) \times \mathbb{L}^2(\Omega)$) to (\ref{form1a}) with $\partial \Omega_T = \emptyset$. 
%%%%%
\subsection{Fine-scale discretization and Picard iteration for linearization}
\label{finedis_e}

The solution $\bfa{u} \in \bfa{V}$ to (\ref{form3}) satisfies
\begin{equation}\label{2cem}
 q(\bfa{u},\bfa{v}) = (\bfa{f},\bfa{v}), \quad \forall \bfa{v} \in \bfa{V}\,,
\end{equation}
where
\begin{align}\label{3cem}
 q(\bfa{u}, \bfa{v}) = \int_{\Omega} 
 \bfa{a} (\bfa{x}, \bfa{Du}) \cdot \bfa{Dv} \, \dx, 
 \quad (\bfa{f},\bfa{v}) = \int_{\Omega} \bfa{f} \cdot \bfa{v} \, \dx\,.
\end{align}
Here, $(\cdot , \cdot)$ represents the standard inner product.

Starting with an initial guess $\bfa{u}^0 = \bfa{0}$, to solve the equation (\ref{form3}), we will linearize it by the Picard iteration, that is, we solve
\begin{align}
 -\ddiv (\kappa(\bfa{x}, |\bfa{D}( \bfa{u}^n)|) 
 \bfa{D}( \bfa{u}^{n+1})) &= \bfa{f} \quad \text{in } \Omega\,, \label{4cem} \\
 \bfa{u}^{n+1} &= \bfa{0} \quad \text{on } \partial \Omega\,, \label{4cemn1}
\end{align}
where superscripts involving $n\text{ } (\geq 0)$ denote respective iteration levels.  

To discretize (\ref{4cem})-(\ref{4cemn1}), we use the notion of fine grid $\Oh$ and coarse grid $\PH$ as well as their related definitions from Section \ref{pre} and Subsection \ref{overv}, respectively.  

On the fine grid $\mathcal{T}_h$, we will approximate the solution of (\ref{2cem}), denoted by $\bfa{u}_h$ (or $\bfa{u}$ for simplicity).  Toward describing the details of the Picard iteration algorithm, we define the bilinear form 
$a(\cdot, \cdot ; \cdot)$:
\begin{equation}\label{biform}
 a(\bfa{u},\bfa{v}; |\bfa{Dw}|) = \int_{\Omega} \kappa
 (\bfa{x},|\bfa{Dw}|) 
 (\bfa{D} \bfa{u} \cdot \bfa{D} \bfa{v}) \dx
\end{equation}
and the functional $J(\cdot)$:
\begin{equation}\label{func}
 J(\bfa{v}) = \int_{\Omega} \bfa{f} \cdot \bfa{v} \dx\,.
\end{equation}

Given $\bfa{u}_h^n$, 
the next approximation $\bfa{u}_h^{n+1}$ is 
the solution of the linear elliptic equation 
\begin{equation}\label{leq1}
 a(\bfa{u}_h^{n+1},\bfa{v}; 
 |\bfa{D}(\bfa{u}_h^n)|) = J(\bfa{v}), \quad \forall \bfa{v} 
 \in \bfa{V}_h\,.
\end{equation}
This is an approximation of the linear equation 
\begin{equation}\label{leq2}
 -\ddiv (\kappa(\bfa{x},|\bfa{D}( \bfa{u}_h^n)|) 
 \bfa{D}( \bfa{u}_h^{n+1})) = \bfa{f}\,.
\end{equation}   

We reformulate the iteration (\ref{leq1}) in a matrix form.  
%When the bilinear form $a(\cdot, %\cdot ; |\bfa{D}(\bfa{u}^n)|)$ is symmetric 
%and positive-definite, 
That is, we define $\bfa{A}_h^n$ by 
\begin{equation}\label{matrix1}
 a(\bfa{w},\bfa{v};|\bfa{D}( \bfa{u}_h^n)|) = 
 \bfa{v}^{\text{T}}\bfa{A}_h^n \bfa{w} \quad \forall \bfa{v},\bfa{w} 
 \in \bfa{V}_h\,.
\end{equation}
and define vector $\bfa{b}_h$ by 
\begin{equation}\label{vector}
 J(\bfa{v}) = \bfa{v}^{\text{T}} \bfa{b}_h, \quad \forall \bfa{v} 
 \in \bfa{V}_h\,.
\end{equation}
In particular, let
\begin{equation}\label{basisVh}
 \{\bfa{p}_1, \cdots, \bfa{p}_c\}
\end{equation}
 be an orthonormal basis for $\bfa{V}_h$.  Then, $\bfa{b}_h$ is exactly the vector whose the $i$th component is $(\bfa{f}, \bfa{q}_i)$, and $\bfa{A}^n_h$ is a symmetric, positive definite matrix with 
 \begin{equation}\label{An}
  \bfa{A}^n_{ij,h} = a(\bfa{p}_j, \bfa{p}_i; |\bfa{Du}^n_h|)\,.
 \end{equation} 
%for b = b_1 q_1 + ... + b_i q_i + ... + b_c q_c, use orthonormality of q_i, (q_i,b) = \int q_i \cdot b = (q_i, b_i q_i) = b_i
%$a_n(u^{n+1},q_i) = (f,q_i)=q_i^T A^n u^{n+1} = q_i^T b = (F)_i = (f, q_i)$.
Thus, in $\bfa{V}_h$, Eq.\ (\ref{leq1}) can be rewritten in the following matrix form:
\begin{equation}\label{matrix2}
 \bfa{A}^n_h \bfa{u}_h^{n+1} = \bfa{b}_h\,.
\end{equation}

Furthermore, at the $(n+1)$th Picard iteration, we can solve Eq.\ (\ref{matrix2}) for the multiscale solution $\bfa{u}^{n+1}_{\tu{ms}} \in \bfa{V}_{\tu{ms}}$ by using the CEM-GMsFEM (to be discussed in the next Subsections \ref{cemoff} and \ref{onba}), with multiscale basis functions for $\bfa{V}_{\tu{ms}}$ computed 
%for $|\bfa{D} \bfa{u}_{\tu{ms}}^n|$ 
in each coarse region $w_i, i = 1, \cdots, N_v$.  
%This means that we will construct multiscale space $\bfa{V}_{ms} (\subset \bfa{V})$, where the solution is obtained.  
%That is, we seek $\bfa{u}_{ms} \in \bfa{V}_{ms}$ such that
%\begin{equation}\label{5cem}
% a(\bfa{u}_{ms},\bfa{v}) = (\bfa{f},\bfa{v}), 
% \quad \forall \bfa{v} \in \bfa{V}_{ms}\,.
%\end{equation}

Each of $\bfa{u}_h$ and $\bfa{u}_{\tu{ms}}$ is computed in a separate Picard iteration procedure, whose termination criterion is that the relative $\bfa{L}^2$ difference is less than $\delta_0$, which can be found in Subsection \ref{onba} and Section \ref{nume} ($\delta_0=10^{-5}$).
%relative criteria

%%%%
\subsection{CEM-GMsFEM for nonlinear elasticity problem}\label{cemforn}
\subsubsection{Overview}\label{overve}
We will construct the offline and online spaces.  As in \cite{gne}, we will focus on the effects of the nonlinearities.  From the linearized equation (\ref{leq2}), we can define offline multiscale basis functions (following the framework of the CEM-GMsFEM) and construct online multiscale basis functions (based on an adaptive enrichment algorithm).

Given $\bfa{u}^n$ (which can represent either $\bfa{u}_h^n$ or $\bfa{u}^n_{\tu{ms}}$, context-dependently).  At the considering $(n+1)$th Picard iteration, we will get the fine-scale solution $\bfa{u}_h^{n+1} \in \bfa{V}_h$ by solving the variational problem
\begin{equation}\label{f3}
 a_n(\bfa{u}_h^{n+1}, \bfa{v}) = (\bfa{f}, \bfa{v}), \quad 
\forall \bfa{v} \in \bfa{V}_h\,,
\end{equation}
where
\begin{equation}\label{biform1e}
 a_n(\bfa{w},\bfa{v}) = \int_{\Omega} \kappa
 (\bfa{x},|\bfa{Du}^n|) 
 (\bfa{D} \bfa{w} \cdot \bfa{D} \bfa{v}) \dx\,.
\end{equation} 

At the $n$th Picard iteration, the space $\bfa{V}_h$ is equipped with the energy norm $\|\bfa{v}\|^2_{\bfa{V}_h} = a_n(\bfa{v}, \bfa{v})$.

%%%%
\subsubsection{General idea of the CEM-GMsFEM for nonlinear elasticity problem}\label{gideae}
The general idea here is as in the dynamic case (Subsection \ref{gidea}).  In this static case, at the current $n$th Picard iteration, we will use the continuous Galerkin (CG) formulation, with a similar form to the fine-scale problem (\ref{f3}).  More specifically, at the $m$th inner iteration, we will construct the multiscale space $\bfa{V}^m_{\tu{ms}} (\subset \bfa{V})$.  That is, we seek $\bfa{u}^m_{\tu{ms}} \in \bfa{V}^m_{\tu{ms}}$ such that
\begin{equation}\label{5cem}
 a_n(\bfa{u}^m_{\tu{ms}},\bfa{v}) = (\bfa{f},\bfa{v}), 
 \quad \forall \bfa{v} \in \bfa{V}^m_{\tu{ms}}\,.
\end{equation}

%in computation, we will solve on fine-grid.
We remark that $\bfa{u}^m_{\tu{ms}}$ from the above problem is in a continuous space.  In numerical simulations, at the current $n$th Picard iteration, we will use the first-order finite elements on the fine grid 
$\mathcal{T}_h$ to compute the multiscale basis functions.  Each multiscale basis function then can be treated as a column vector 
$\bfa{\Phi}_i$.  Let $\bfa{P}=[\bfa{\Phi}_1, \cdots, \bfa{\Phi}_{Nms}]$ be the matrix that is formed by all $Nms$ multiscale basis functions (at the $m$th inner iteration).  Hence, the multiscale solution satisfies $\bfa{u}^m_{\tu{ms}} = (\bfa{P}^{\tu{T}} \bfa{A}^n_h \bfa{P})^{-1} (\bfa{P}^{\tu{T}} \bfa{b}_h)$ in $\bfa{V}^m_{\tu{ms}}$.  Projecting the coarse solution $\bfa{u}^m_{\tu{ms}}$ onto $\bfa{V}_h$, we obtain $\bfa{u}^f_{\tu{ms}} = \bfa{P}\bfa{u}^m_{\tu{ms}}$.

Our results show that the combination of offline and online multiscale basis functions (via adaptive enrichment) within the CEM-GMsFEM will give a faster convergence of the sequence of multiscale solutions $\{\bfa{u}^m_{\tu{ms}} \}_{m\geq 1}$ to the fine-scale solution $\bfa{u}_h$ than within the GMsFEM in \cite{gne}.

%%%%%
\subsection{Construction of CEM-GMsFEM offline multiscale basis functions}\label{cemoff}
The readers who have already gone through Sections \ref{cembase} for the dynamic case may skip this Subsection \ref{cemoff}, which are similar to Subsections \ref{auxb} and \ref{mss}.

Toward clarity for the static case, we still present here this Subsection \ref{cemoff} regarding the construction of the offline multiscale basis functions, at the $n$th Picard iteration ($n \geq 0$).  There are two stages.  The first stage is to construct the auxiliary multiscale basis functions in the framework of the GMsFEM.  The second stage is to construct the offline multiscale basis functions 
by solving some constraint energy minimizing (CEM) problems in the oversampled region.  

%%%%
\subsubsection{Auxiliary multiscale basis functions}\label{auxbe}
In each coarse block $K_i$, the auxiliary multiscale basis functions are constructed by solving a spectral problem.  More specifically, for each coarse block $K_i$, we let $\bfa{V}(K_i)$ be the restriction of $\bfa{V}$ on $K_i$.  Then, we solve the local spectral problem:  
find $(\lambda^i_j, \bfa{\phi}^i_j) \in \mathbb{R} \times \bfa{V}(K_i)$ ($j=1,2,\cdots$) such that
\begin{equation}\label{lsp}
a^i_n(\bfa{\phi}^i_j, \bfa{w}) = \lambda^i_j s^i_n(\bfa{\phi}^i_j,\bfa{w})\,, \quad 
\forall \bfa{w} \in \bfa{V}(K_i)\,,
\end{equation}
where 
\begin{equation}\label{ai}
 a^i_n(\bfa{v}, \bfa{w}) = 
\int_{K_i} \kappa(\bfa{x}, |\bfa{Du}^n_{\tu{ms}}|)\bfa{Dv} \cdot \bfa{Dw} \, \dx\,,
\end{equation}
%in Picard iteration, we start at initial guess u^{old}_ms = u^n_ms and compute V^{old}_ms only one time as we do not update the basis after each Picard step
and
\begin{equation}\label{si}
s^i_n(\bfa{v}, \bfa{w}) = 
\int_{K_i} \tilde{\kappa} \bfa{v} \cdot \bfa{w}\, \dx\,, 
\end{equation}
in which,
\[\tilde{\kappa}=  
\kappa(\bfa{x}, |\bfa{Du}^n_{\tu{ms}}|)\sum_{k=1}^{N_v}|\nabla  \chi_k|^2\,,\]
and $\{\chi_k\}$ is a set of partition of unity functions (see \cite{unitybabu}) with respect to the coarse grid.  
%One can take $\{\bfa{\chi}_k\}$ as the standard multiscale basis functions or the standard piecewise linear functions.
Our hypothesis is that the eigenfunctions satisfy the normalized condition 
\[s^i_n(\bfa{\phi}_j^i, \bfa{\phi}_j^i)=1\,.\]  We still denote by $\lambda^i_j$ the eigenvalues of (\ref{lsp}) arranged in nondecreasing order.  
Then, using the first $L_i$ corresponding eigenfunctions, we will construct our local auxiliary multiscale space $\bfa{V}^i_{\tu{aux}}$, where
\[\bfa{V}_{\tu{aux}}^i = \textup{span}\{\bfa{\phi}^i_j: 1 \leq j \leq L_i\}\,.\]
Also, let $\Lambda$ be the minimum of the first discarded eigenvalues, that is 
\begin{equation}\label{Lam}
 \Lambda = \min_{1 \leq i \leq N} 
 \lambda^i_{L_i +1}\,,
\end{equation}
where $\lambda^i_{L_i +1} = \mathcal{O}(1)$ by construction (\cite{cem2f}).  In global setting, the auxiliary space $\bfa{V}_{\tu{aux}}$ is determined by the sum of all local auxiliary spaces 
$\bfa{V}_{\tu{aux}}^i= \bfa{V}_{\tu{aux}}(K_i)$:
\[\bfa{V}_{\tu{aux}} = \bigoplus_{i=1}^N \bfa{V}_{\tu{aux}}^i\,.\]

Given a local auxiliary multiscale space $\bfa{V}^i_{\tu{aux}}$, the bilinear form 
$s_n^i$ in (\ref{si}) leads to an inner product with norm
\[\|\bfa{v}\|_{s_n^i} = \sqrt{s_n^i(\bfa{v}, \bfa{v})} \,.\]
We thus define
\[s_n(\bfa{v},\bfa{w}) = \sum_{i=1}^N s_n^i(\bfa{v},\bfa{w})\,,
\quad \|\bfa{v}\|_{s_n}=\sqrt{s_n(\bfa{v},\bfa{v})}\,, \quad \forall 
\bfa{v} \in \bfa{V}_{\tu{aux}}\,.\]

In the continuous space $\bfa{V}$, given a function $\bfa{\phi}_j^i \in \bfa{V}_{\tu{aux}}$, we introduce the notion of 
$\bfa{\phi}_j^i$-orthogonality: a function $\bfa{\psi} \in \bfa{V}$ is called
$\bfa{\phi}_j^i$-orthogonal if 
\[s_n(\bfa{\psi},\bfa{\phi}_j^i)=1\,, \quad 
s_n(\bfa{\psi},\bfa{\phi}_{j'}^{i'})=0 \quad \text{ if } j'\neq j 
\text{ or } i' \neq i\,.\]
We now define a projection operator $\pi^i_n$ from space $\bfa{V}(K_i)$ 
%as V^i \subset L^2(K_i)
to $\bfa{V}^i_{\tu{aux}}$ as follows: 
\[\pi^i_n(\bfa{u}) = \sum_{j=1}^{L_i}s_n^i(\bfa{u},\bfa{\phi}_j^i) 
\bfa{\phi}_j^i\,, \quad \forall \bfa{u} \in \bfa{V}(K_i)\,.\]
Furthermore, we let $\pi_n: \bfa{V} \to \bfa{V}_{\tu{aux}}$ be the 
projection with respect to the inner product $s_n(\bfa{u}, \bfa{w})$.  Then, 
we define the operator $\pi_n$ by
\[\pi_n(\bfa{u}) = \sum_{i=1}^N \sum_{j=1}^{L_i} 
s_n^i(\bfa{u},\bfa{\phi}_j^i) \bfa{\phi}_j^i\,, \quad
\forall \bfa{u} \in \bfa{V}\,.\]
Note that $\pi_n = \dd \sum_{i=1}^N \pi^i_n$.  The kernel of the operator $\pi_n$ 
restricted to $\bfa{V}$ is denoted by 
\[\tilde{\bfa{V}} = 
\{\bfa{w} \in \bfa{V} | \pi_n(\bfa{w}) = 0\}\,.\]

%%%%%
\subsubsection{Offline multiscale basis functions}\label{offb}
After building the auxiliary space, we can construct offline multiscale basis functions for the iteration $n$ ($\geq 0$).  Given a coarse 
block $K_i$, we define an oversampled domain $K_{i,k} \subset \Omega$ by expanding $K_i$ by $k$ coarse-grid layers ($k \geq 1$ is an integer).  For each $\bfa{\phi}^i_{j} \in \bfa{V}_{\tu{aux}}$, we define the multiscale basis function 
$\bfa{\psi}^{i,\tu{ms}}_{j} \in \bfa{V}(K_{i,k})$ by
\begin{equation}\label{msbf1}
 \bfa{\psi}^{i,\tu{ms}}_{j} = \textup{argmin} 
 \{a_n(\bfa{\psi},\bfa{\psi}) \, | \, \bfa{\psi} 
 \in \bfa{V}(K_{i,k})\,, \bfa{\psi} 
 \text{ is } \bfa{\phi}^i_j \textup{-orthogonal} \}\,,
\end{equation}
where $\bfa{V}(K_{i,k}) = \bfa{H}_0^1(K_{i,k})$.  Using Lagrange Multiplier, we can rewrite the problem 
(\ref{msbf1}) as follows:  find $\bfa{\psi}^{i,\tu{ms}}_j \in 
\bfa{V}(K_{i,k})$ and $\bfa{\nu} \in \bfa{V}^i_{\tu{aux}}$ such that
\begin{align}\label{lm}
 \begin{split}
  a_n(\bfa{\psi}^{i,\tu{ms}}_j,\bfa{p}) + s_n(\bfa{p},\bfa{\nu}) &= 0 \quad \forall \bfa{p} \in \bfa{V}(K_{i,k})\,,\\ 
  s_n(\bfa{\psi}^{i,\tu{ms}}_j - \bfa{\phi}^i_j,\bfa{q}) &= 0 \quad \forall \bfa{q} \in \bfa{V}_{\tu{aux}}(K_{i,k})\,,
 \end{split}
\end{align}
where $\bfa{V}_{\tu{aux}}(K_{i,k})$ is the union of all local auxiliary spaces for $K_r \subset K_{i,k}$.    
%see the equivalence in p10, step 1, in \cite{cem1}

This continuous problem can be solved numerically within the fine-scale mesh $\bfa{V}_h$, at the current $n$th Picard iteration.  In particular, let $\bfa{M}^n_h$  be the matrix such that $\bfa{M}^n_{ij,h} = s_n(\bfa{p}_j,\bfa{p}_i)\,,$ where $\bfa{p}_j,\bfa{p}_i$ are from (\ref{basisVh}).  
Restricting $\bfa{A}^n_h$ from (\ref{An}) and the above $\bfa{M}^n_h$ on $K_{i,k}$ , we respectively obtain $\bfa{A}^i_h$ and $\bfa{M}^i_h\,,$ where the superscript $n$ is omitted.  Then, let $\bfa{P}^i$ be the matrix that consists of all the 
discrete auxiliary basis functions for the space $\bfa{V}_{\tu{aux}}(K_{i,k})$.

The problem (\ref{lm}) can be recast as the 
following matrix
\begin{equation}\label{lcm}
\begin{pmatrix}
\bfa{A}^i_h& \bfa{M}^i_h\bfa{P}^i\\
(\bfa{M}^i_h\bfa{P}^i)^{\tu{T}}& \bfa{0}
\end{pmatrix}  
\begin{pmatrix}
\bfa{\psi}^i_h\\
\bfa{\nu}^i_h
\end{pmatrix} = 
\begin{pmatrix}
\bfa{0}\\
\bfa{I}_i
\end{pmatrix} \,,
\end{equation}
where $\bfa{P}^i_j$ is the $j$th column of $\bfa{P}^i$, $\bfa{\psi}^{i}_{j,h}$ is the discretization of $\bfa{\psi}^{i,\tu{ms}}_j$, $\bfa{I}_i$ is a sparse matrix whose nonzero elements (all are 1) are in the diagonal of the matrix, and the nonzero elements' positions depend on the index order of $K_i$ in $K_{i,k}$ (\cite{cemgle}).

Thanks to \cite{cem1}, for each $\bfa{\phi}^i_{j} \in \bfa{V}_{\tu{aux}}$, from the $\bfa{\phi}^i_j-$orthogonality in (\ref{msbf1}), we obtain a relaxed version of the multiscale basis functions.  That is, we solve the following un-constrainted minimization problem:  find multiscale basis function $\bfa{\psi}^{i,\tu{ms}}_j \in \bfa{V}(K_{i,k})$ such that
\begin{equation}\label{uc}
 \bfa{\psi}^{i,\tu{ms}}_j = \textup{argmin} 
 \{a_n(\bfa{\psi}, \bfa{\psi}) + s_n(\pi_n(\bfa{\psi}) - 
 \bfa{\phi}^i_j, \pi_n(\bfa{\psi}) - \bfa{\phi}^i_j) \, 
 | \, \bfa{\psi} \in \bfa{V}(K_{i,k}) \}\,,
\end{equation}
which is equivalent to the following variational formulation
\begin{equation}\label{vuc}
 a_n(\bfa{\psi}^{i,\tu{ms}}_j,\bfa{v}) + s_n(\pi_n(\bfa{\psi}^{i,\tu{ms}}_j),\pi_n(\bfa{v})) = s_n(\bfa{\phi}^i_j,\pi_n(\bfa{v}))\,, \quad \forall \bfa{v} \in \bfa{V}(K_{i,k})\,.
\end{equation}
With the same notation as above, Eq.\ (\ref{vuc}) has the following matrix formulation:
\begin{equation}\label{mvuc}
 (\bfa{A}^i_h + \bfa{M}^i_h(\bfa{P}^i\, \bfa{P}^{i,\tu{T}})\bfa{M}_h^{i,\tu{T}}) \, \bfa{\psi}^{i}_{j,h} = \bfa{P}^i_j \bfa{M}_h^{i,\tu{T}}\,.
\end{equation}

For each auxiliary multiscale basis function $\bfa{\phi}^i_j \in \bfa{V}_{\tu{aux}}$, one can obtain a multiscale basis function $\bfa{\psi}^{i,\tu{ms}}_j$.  
Finally, the span of these multiscale basis functions forms the multiscale finite element space
\[\bfa{V}_{\tu{ms}} := \tu{span} \{ \bfa{\psi}^{i,\tu{ms}}_j: 1 \leq j \leq L_i, 1 \leq i \leq N \}\,.\]  
This method is thus called CEM-GMsFEM because the construction of 
the multiscale basis includes solving spectral problems and 
energy minimization problems.  The (local) multiscale basis functions $\bfa{\psi}^{i,\tu{ms}}_j \in \bfa{V}(K_{i,k})$ are used to approximate the related global multiscale basis functions 
$\bfa{\psi}^i_j \in \bfa{V}$, which is defined in the same manner (\cite{cem1}), 
that is to say,
\begin{equation}\label{gcmbf}
 \bfa{\psi}^i_j = \textup{argmin} \{a_n(\bfa{\psi}, \bfa{\psi}) \, | \, \bfa{\psi} \in \bfa{V}\,, \bfa{\psi} \text{ is } \bfa{\phi}^i_j
 \textup{-orthogonal}\}
\end{equation}
for the constraint case, and 
\begin{equation}\label{grmbf}
 \bfa{\psi}^i_j = \textup{argmin} \{a_n(\bfa{\psi}, \bfa{\psi}) + s_n(\pi_n(\bfa{\psi}) - \bfa{\phi}^i_j, 
 \pi_n(\bfa{\psi}) - \bfa{\phi}^i_j) \, | \, \bfa{\psi} \in \bfa{V}\}
\end{equation}
for the relaxed case, which is equivalent to the following global problem (see \cite{cem1}):  find $\bfa{\psi}^i_j \in \bfa{V}$ such that
\[a_n(\bfa{\psi}^i_j,\bfa{v}) + s_n(\pi_n(\bfa{\psi}^i_j),\pi_n(\bfa{v})) = s_n(\bfa{\phi}^i_j, \pi_n(\bfa{v}))\,, \quad \forall \bfa{v} \in \bfa{V}\,.\]
The global multiscale finite element space is now defined by 
\[\bfa{V}_{\tu{glo}} = \textup{span} \{\bfa{\psi}^i_j \, | \, 1 \leq j \leq L_i\,, 1 \leq i \leq N \}\,.\]
These global basis functions have an exponential decay 
property (\cite{cem1}), 
 which motivates the definitions of the multiscale basis functions $\bfa{\psi}^{i, \tu{ms}}_j$ (\ref{msbf1}) having local supports (\cite{cem2f}).
 %, which is important to the convergence analysis of the following online adaptive enrichment method.

%%%%%
\subsection{Online multiscale basis functions and adaptive enrichment}\label{onba}
Now, we will introduce an online enrichment process for this CEM-GMsFEM, at the $n$th Picard iteration.  First, 
the construction of online multiscale basis functions is shown.  
Second, an adaptive enrichment method based on an error estimate is 
presented.  

The online basis functions, in online stage, are constructed iteratively 
using the residual of previous multiscale solution, which contains 
the source and global information of the media.
%We remark that the error will decay rapidly such that the error will 
%be within an acceptable range in the first or two iterations.

At the current $n$th Picard iteration, we are given a coarse neighborhood $w_i$, an inner adaptive iteration $m$th, and an approximation space $\bfa{V}^m_{\tu{ms}}$.  Recall that the GMsFEM solution $\bfa{u}^m_{\tu{ms}} \in \bfa{V}^m_{\tu{ms}}$ ($\subset \bfa{V}$) can be obtained by solving (\ref{5cem}):
\begin{align*}
 a_n(\bfa{u}^m_{\tu{ms}}, \bfa{v}) = (\bfa{f}, \bfa{v})\,, 
 \quad \forall \bfa{v} \in \bfa{V}^m_{\tu{ms}}\,.
\end{align*}  
A residual functional $r: \bfa{V} \to \mathbb{R}$ is then defined by 
\begin{equation}\label{res}
 r(\bfa{v}) = a_n(\bfa{u}^m_{\tu{ms}},\bfa{v}) - 
 \int_{\Omega} \bfa{f} \cdot \bfa{v} \,, 
 \quad \forall \bfa{v} \in \bfa{V}\,,
\end{equation}
whose discretization in matrix form is
%A_h u_h - A_h u^f_{ms}, u^f_{ms} = R u_{ms}, F_h = A_h u_h
%u_{ms} = (R^T A_h R)^{-1} (R^T F_h)
\[\bfa{b}_h -\bfa{A}^n_h(\bfa{P}((\bfa{P}^T \bfa{A}^n_h \bfa{P})^{-1}(\bfa{P}^{\tu{T}}\bfa{b}_h)))\,.\]
Given a coarse neighborhood $w_i$, for all $\bfa{v} \in \bfa{V}$, 
we define the local residual functional 
$r_i: \bfa{V} \to \mathbb{R}$ by
\[r_i(\bfa{v}) = r(\chi_i \bfa{v})\,,\]
which gives a measure of the error $\bfa{u} - \bfa{u}^m_{\tu{ms}}$ 
in $w_i$.

Let $w_i^{+}$ be an extending of $w_i$ by a few coarse blocks.  
Using the local residual $r_i$, we can construct 
online basis function $\bfa{\beta}^i_{\tu{ms}}$ whose 
support is an oversampled region $w_i^{+}$.  In particular, 
the online basis function $\bfa{\beta}^i_{\tu{ms}} \in \bfa{H}^1_0(w_i^{+})$ is 
obtained by solving the following equation:
\begin{equation}\label{re1}
 a_n(\bfa{\beta}^i_{\tu{ms}},\bfa{v}) + s_n(\pi_n(\bfa{\beta}^i_{\tu{ms}}),\pi_n(\bfa{v})) 
 = r_i(\bfa{v})\,, \quad \forall \bfa{v} \in \bfa{H}^1_0(w_i^{+})\,.
\end{equation}
Solving Eq.\ (\ref{re1}) is similar to solving Eq.\ (\ref{vuc}).  
The online multiscale basis function is also a localization result of 
the corresponding global online basis function $\bfa{\beta}^i_{\tu{glo}} \in 
\bfa{V}$ defined by 
\begin{equation}\label{re2}
 a_n(\bfa{\beta}^i_{\tu{\tu{glo}}} , \bfa{v}) + s_n(\pi_n(\bfa{\beta}^i_{\tu{glo}}),\pi_n(\bfa{v}))
 =r_i(\bfa{v})\,, \quad \forall \bfa{v} \in \bfa{V}\,.
\end{equation}
In practice, one can perform the above construction 
%for some 
%selected $r_i$ (with $i \in I$ for an index set $I$) 
based on an adaptive 
criterion.  After constructing the online basis functions, we can enrich 
the offline multiscale space by adding the online basis: 
\[\dd \bfa{V}^{m+1}_{\tu{ms}} = \bfa{V}^m_{\tu{ms}} + \textup{span}_{1\leq i \leq N_v} 
\{\bfa{\beta}^i_{ms}\}\,.\]
Within this new multiscale finite element space, we can compute new 
multiscale solution by solving Eq.\ (\ref{5cem}).
%and compute 
%new multiscale finite element space.  
%This process can be repeated to enrich our multiscale 
%space until the residual norm is smaller than a given tolerance.  
Before presenting the online adaptive enrichment algorithm, we first 
define the $a_n$-norm $\| \cdot \|_{a_n}$, where $\| \bfa{w} \|_{a_n}^2 
= a_n(\bfa{w}, \bfa{w})$.

%%%%%

\subsubsection{Online adaptive enrichment algorithm}\label{oaa}
Assume that we are at the $n$th Picard step.  First, we choose an initial space when the inner iteration $m=0$, that is, $\bfa{V}^0_{\tu{ms}}$, which is obtained by using the offline multiscale basis functions constructed in Subsection \ref{cemoff}. 

For each inner iteration $m=0,1, \cdots,$ we assume that $\bfa{V}^m_{\tu{ms}}$ is given.  Then, the following procedure allows us to find the new multiscale finite element space $\bfa{V}^{m+1}_{\tu{ms}}$. 

\noindent \textbf{Step 1}:  Find the multiscale solution in the current space $\bfa{V}^m_{\tu{ms}}$.  That is, find $\bfa{u}^m_{\tu{ms}} \in \bfa{V}^m_{\tu{ms}}$ such that
\begin{equation}\label{sin}
 a_n(\bfa{u}^m_{\tu{ms}},\bfa{v}) = (\bfa{f}, \bfa{v})\,, \quad \forall \bfa{v} \in \bfa{V}^m_{\tu{ms}}\,.
\end{equation}

\noindent \textbf{Step 2}:  Construct the local online basis functions.  For each $1 \leq i \leq N_v$ and coarse neighborhood $w_i$, we find online basis function $\bfa{\beta}^i_{\tu{ms}} \in \bfa{H}^1_0(w_i^{+})$ satisfying
\[a_n(\bfa{\beta}^i_{\tu{ms}},\bfa{v}) + s_n(\pi_n (\bfa{\beta}^i_{\tu{ms}}), \pi_n(\bfa{v})) = r^m_i(\bfa{v})\,, \quad \forall \bfa{v} \in \bfa{H}^1_0(w_i^{+})\,,\]
where $\dd r_i^m(\bfa{v}) = a_n(\bfa{u}^m_{\tu{ms}}, \chi_i \bfa{v}) - 
\int_{\Omega} \bfa{f} \cdot (\chi_i \bfa{v}) \,.$

\noindent \textbf{Step 3}:  Enrich the multiscale finite element space by
\[\bfa{V}_{\tu{ms}}^{m+1} = \bfa{V}_{\tu{ms}}^m + 
\textup{span}_{1 \leq i \leq N_v} \{ \bfa{\beta}^i_{\tu{ms}} \}\,.\]

\noindent \textbf{Step 4:}  If the dimension of $\bfa{V}^{m+1}_{\tu{ms}}$ is as large as desired, then stop.  Otherwise, set $m \gets m+1$ and go back to \textbf{Step 1}.

\bigskip

For Picard iteration procedure, in the numerical Section \ref{nume}, the multiscale finite element space $\bfa{V}^{(n+1)}_{\tu{ms}}$ is not needed to be updated at every Picard iteration step $(n+1)$th.  We choose the initial basis function space when $n=0$, that is, $\bfa{V}^{(0)}_{\tu{ms}}$ (obtained from the Online adaptive enrichment algorithm \ref{oaa} with $n=0$) for all Picard iteration steps.  Our obtained numerical results are already good with this initial basis.  Whereas, updating basis at every Picard iteration is not cheap.

%%%%%

\subsubsection{CEM-GMsFEM for nonlinear elasticity}\label{cemne}
We sum up the main steps (as in \cite{gne}) of using the CEM-GMsFEM to solve the problem (\ref{form3})-(\ref{D}):
select a Picard iteration
stop tolerance value $\delta_0\in \mathbb{R}_{+}$ (where $\delta_0 = 10^{-5}$ and will be presented in Section \ref{nume}). We also choose a starting guess of $\bfa{u}^{\text{old}}_{\tu{ms}}$, 
%n=0 to remember the level of Picard iteration
and compute $\kappa^{\text{old}}(\boldsymbol{x})=\dfrac{1}{1-\beta(\boldsymbol{x})|(\bfa{Du}^{\text{old}}_{\tu{ms}})|}$ 
and the multiscale space  $\bfa{V}^{\text{old}}_{\tu{ms}}=\bfa{V}^{(0)}_{\tu{ms}}$ (obtained from the Online adaptive enrichment algorithm \ref{oaa} with $n=0$), then we repeat the following steps:

\noindent \textbf{Step 1:}  
%n=n+1
Solve for $\bfa{u}^{\text{new}}_{\tu{ms}} \in \bfa{V}^{\text{old}}_{\tu{ms}}$ from the equation (as (\ref{5cem})) as follows: 
\begin{equation}\label{eold}
 a_{\tu{old}}(\bfa{u}^{\tu{new}}_{\tu{ms}}, \bfa{v}) = (\bfa{f}, \bfa{v}) \quad \forall \bfa{v} \in \bfa{V}^{\text{old}}_{\tu{ms}}\,.
\end{equation}
If $\dfrac{\| \bfa{u}^{\text{new}}_{\tu{ms}}-\bfa{u}^{\text{old}}_{\tu{ms}}\|_{\bfa{V}_h}}{\| \bfa{u}^{\text{old}}_{\tu{ms}}\|_{\bfa{V}_h}}>\delta_0$,
let $\bfa{u}^{\text{old}}_{\tu{ms}}=
\bfa{u}^{\text{new}}_{\tu{ms}}$
and go to \textbf{Step 2}. \\

\noindent \textbf{Step 2:}  Calculate $\kappa^{\text{new}}(\boldsymbol{x})=\dfrac{1}{1-\beta(\boldsymbol{x})|(\bfa{Du}^{\text{new}}_{\tu{ms}})|}$ and let $\kappa^{\text{old}}(\boldsymbol{x})=\kappa^{\text{new}}(\boldsymbol{x})\,.$ 

\noindent Then go to \textbf{Step 1}.

%%%%%%
\section{Numerical results}\label{nume}
\begin{figure}[H]
	\centering
	\subfigure[Model 1]{
		\includegraphics[width=3in]{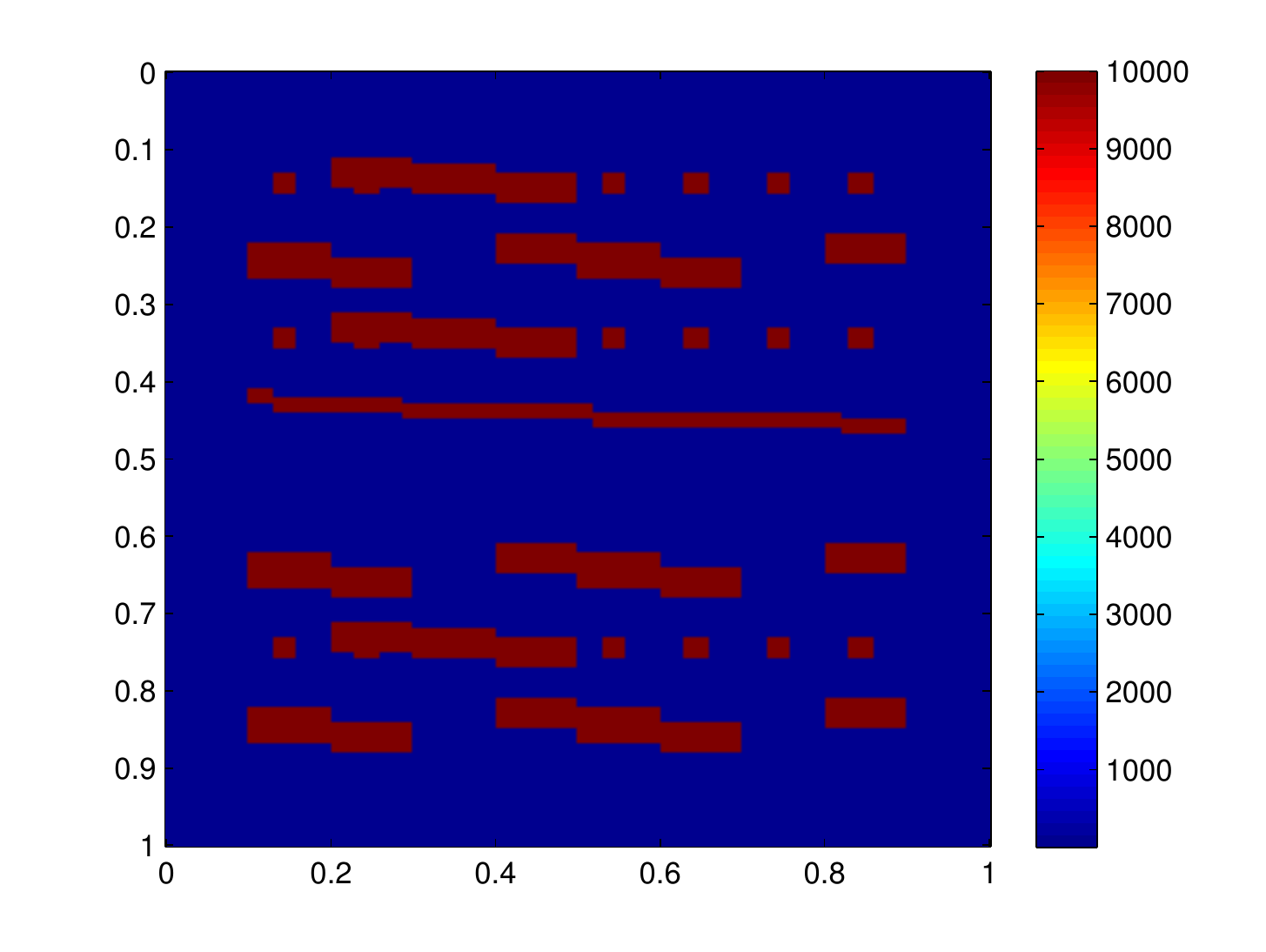}}
	\subfigure[Model 2]{
		\includegraphics[width=3in]{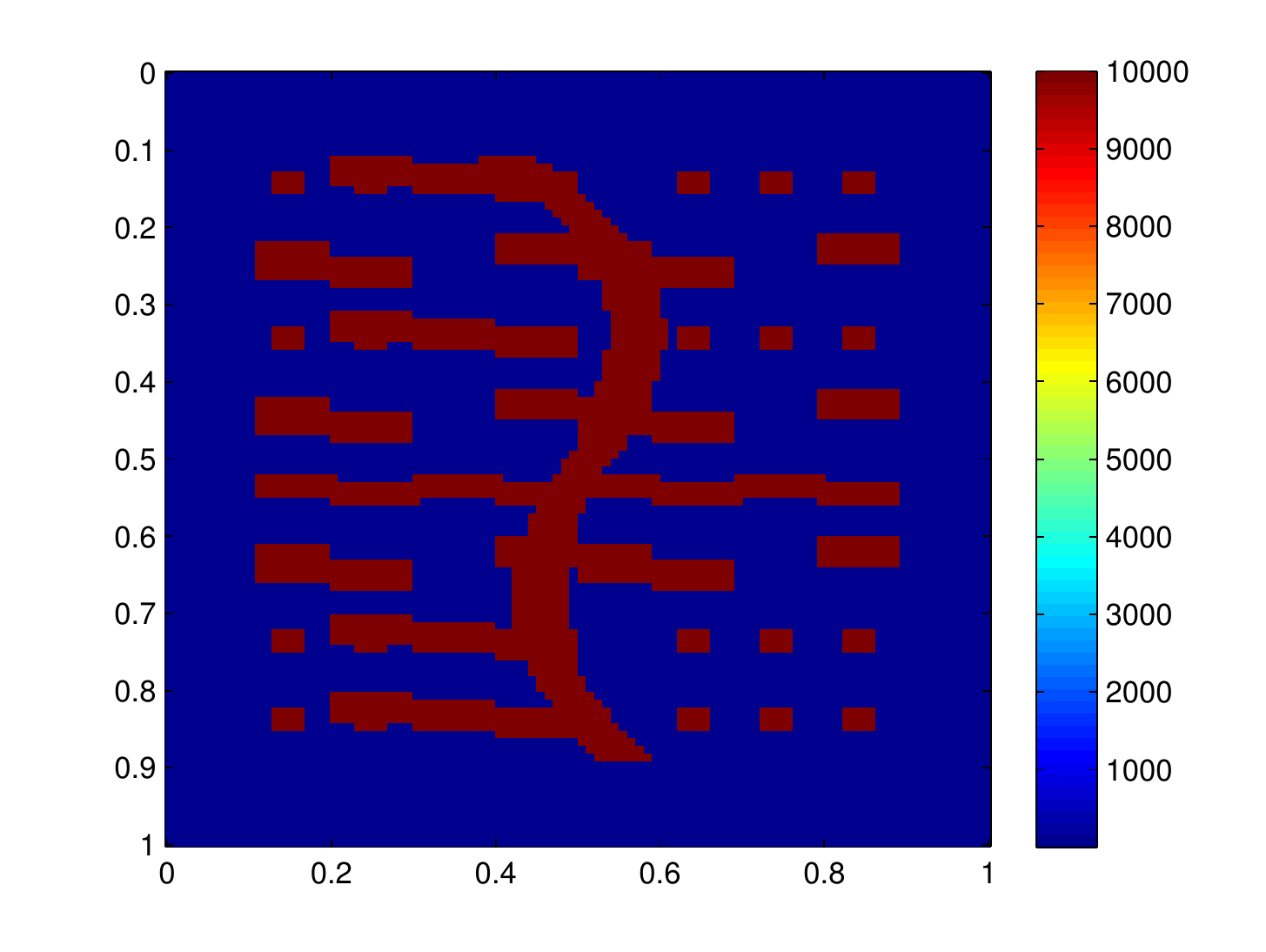}}
		%improved errors reduction
	\caption{Models.}
	\label{fig:model}
\end{figure}

In this section, we will present several numerical experiments to show 
the performance of our method. 
In the simulations, we consider two choices of $\beta(\boldsymbol{x})$,
which are depicted in Figure \ref{fig:model}.
For both test models, the blue region represents $\beta(\boldsymbol{x})=1$
and the red regions represent $\beta(\boldsymbol{x})=10^{4}$. In addition, the precision of the two test models are $200\times 200$, the computational domain
is [0,1]$\times$[0,1]. 
In all tables shown below, $m$ represents the number of oversampling layers,
$J$ is the number of local basis functions, $H$ denotes the coarse-grid size.
If $J=x+y$, then $x$ means the number of offline multiscale basis functions, $y$ represents the 
number of online basis functions.
We take the source term $f=\left(10^{-4}\sqrt{x^2+y^2+1},10^{-4}\sqrt{x^2+y^2+1}\right)$ and $\delta_0 = 10^{-5}$ (for either elasticity or poroelasticity).
\subsection{Static nonlinear elasticity case}\label{nec}
We first consider the static nonlinear elasticity problem.  The CEM-GMsFEM solution will be compared with the fine-grid solution.  At the $(n+1)$th Picard iteration, to quantify the accuracy of our multiscale solutions, 
we use the following
relative weighted $\bfa{L}^2$ error and energy error:
\begin{equation*}
e^u_{\bfa{L}^2}=\frac{||(\bfa{u}_{\tu{ms}}-\bfa{u}_h)||_{\bfa{L}^2(\Omega)}}{||\bfa{u}_h||_{\bfa{L}^2(\Omega)}},\quad
e^u_{a}=\sqrt{\frac{a_n(\bfa{u}_{\tu{ms}}-\bfa{u}_h,\bfa{u}_{\tu{ms}}-\bfa{u}_h)}
	{a_n(\bfa{u}_h,\bfa{u}_h)}}\,, 
\end{equation*}
where the reference solution $\boldsymbol{u}_h$ is computed via (\ref{f3}) on the fine grid, the multiscale solution $\bfa{u}_{\tu{ms}}$ is obtained from (\ref{eold}), and the bilinear form $a_n$ is defined in (\ref{biform1e}).

\begin{figure}[H]
	\centering
	\subfigure{
    \includegraphics{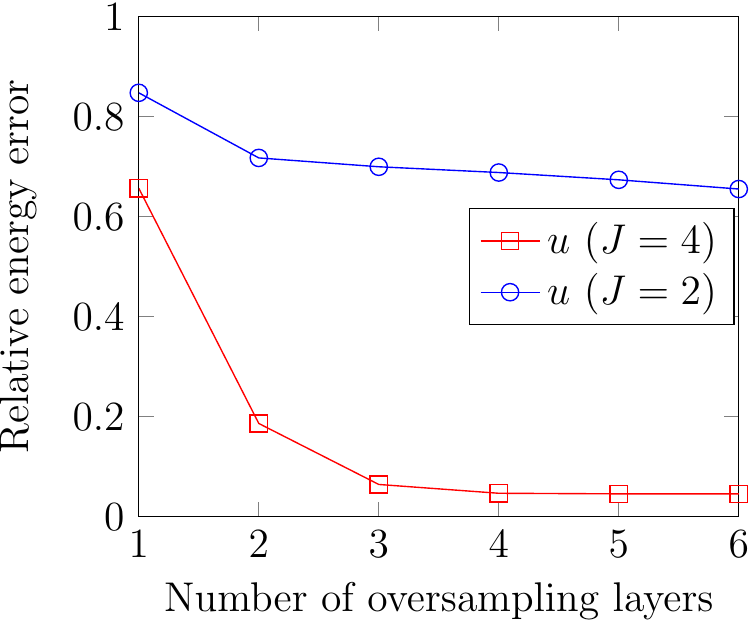}}
    \subfigure{\includegraphics{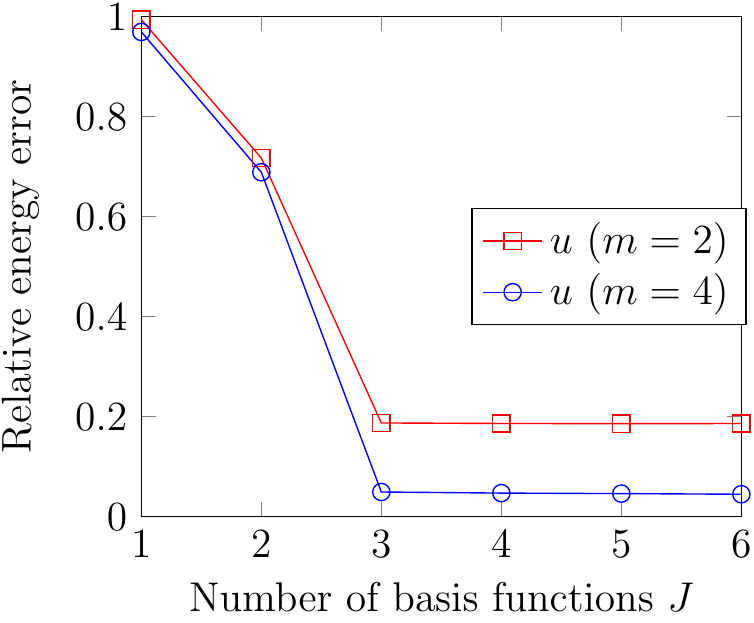}}
	\caption{Relative energy error (Test model 1) for $H=\sqrt{2}/20$ and fixed $J$ (left), fixed $m$ (right).}
	\label{fig:ovnb1}
\end{figure}

\begin{figure}
\centering
	\subfigure{
    \includegraphics{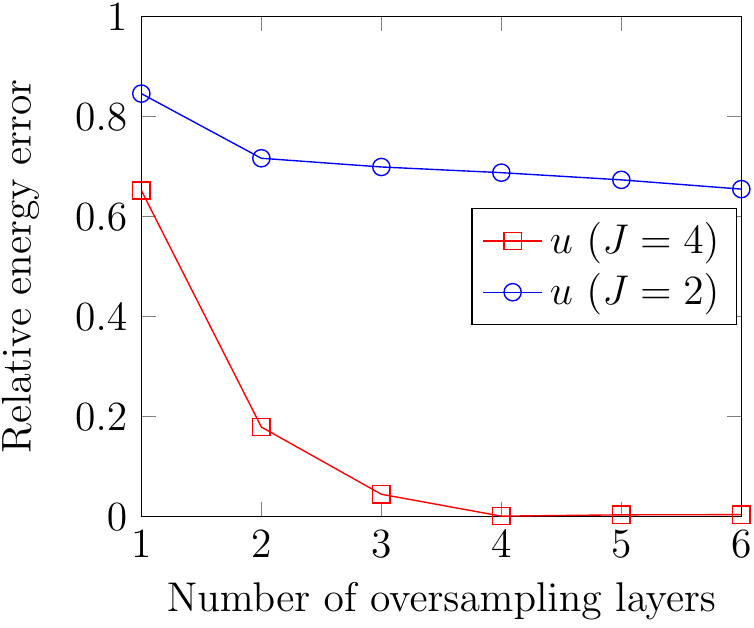}}
    \subfigure{\includegraphics{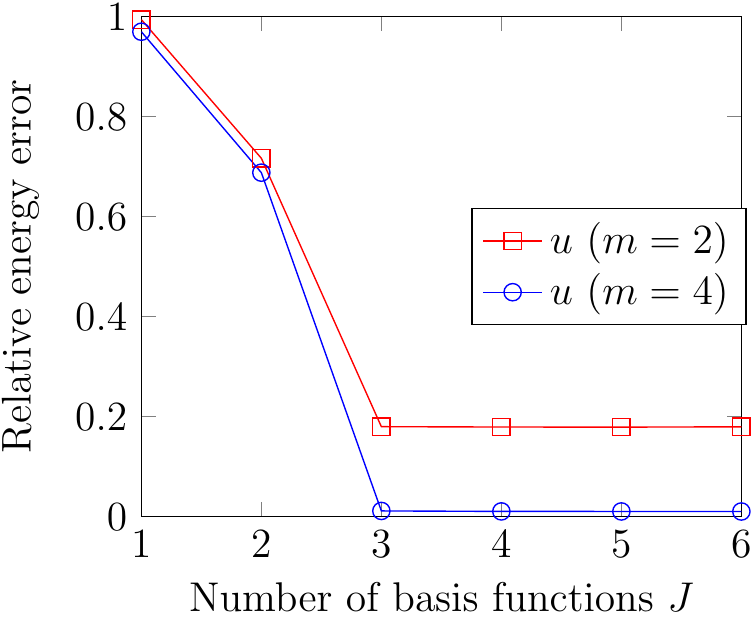}}
	\caption{Relative energy error (Test model 2) for $H=\sqrt{2}/20$ and fixed $J$ (left), fixed $m$ (right).}
	\label{fig:ovnb2}
\end{figure}

First, we study the convergence behavior of the CEM-GMsFEM solution
with respect to the coarse-grid size. We set the number of oversampling layers to $m=3\lfloor\text{log}(H)/\text{log}(\sqrt{2}/10)\rfloor$ and $J=4$ to form the  basis spaces. The results for two test models are shown in Tables \ref{ta:nonlinear_h1} and \ref{ta:nonlinear_h2}, respectively.  We can see clearly for
both test cases that the sequence of CEM-GMsFEM solutions converges as the sequence of coarse-mesh sizes $H$
converges, and it is very accurate.
We also study the effects of oversampling layers and number of basis functions.
The results are plotted in Figure \ref{fig:ovnb1} and Figure \ref{fig:ovnb2}.
It can be observed that increasing the number of basis functions and 
oversampling layers will increase the accuracy of the CEM-GMsFEM solution
as expected. Once $J$ or $m$ exceed some certain numbers, the error 
will no longer decrease.
The performance of using online basis is also investigated,
and the results are presented in Table \ref{ta:nonlinear_online_h1} as well as 
Table \ref{ta:nonlinear_online_h2}. As we can see, the error when 4+2 basis functions are used is less than the error when 6 offline multiscale basis functions are used.  Hence, we can conclude
that residual based online basis functions are more efficient than offline multiscale basis functions.

\begin{table}
	\centering \begin{tabular}{c|c|c||c|c|c}%\hline
		$J$ &$H$& $m$  & $e_{L^2}^u$   & $e_{a}^u$  \tabularnewline\hline
		4&$\sqrt{2}$/10	&3&1.365e-02 & 6.873e-02    \tabularnewline\hline 
		4&$\sqrt{2}$/20	&4&5.864e-03 & 4.680e-02  \tabularnewline\hline
		4&$\sqrt{2}$/40	&5&2.650e-03 & 3.174e-02  \tabularnewline%\hline
	\end{tabular}
	\caption{Numerical results with varying coarse-grid size $H$ for Test model $1$.} 
	\label{ta:nonlinear_h1}
\end{table}

\begin{table}
	\centering \begin{tabular}{c|c|c||c|c|c}%\hline
		$J$ &$H$& $m$  & $e_{L^2}^u$   & $e_{a}^u$  \tabularnewline\hline
		4&$\sqrt{2}$/10	&3&6.986e-04 & 1.405e-02    \tabularnewline\hline 
		4&$\sqrt{2}$/20	&4&2.794e-04 & 1.051e-02 \tabularnewline\hline
		4&$\sqrt{2}$/40	&5&1.319e-04 & 7.564e-03 \tabularnewline%\hline
	\end{tabular}
	\caption{Numerical results with varying coarse-grid size $H$ for Test model $2$.} 
	\label{ta:nonlinear_h2}
	%improved smaller errors
\end{table}
%%%%

\begin{table}
	\centering \begin{tabular}{c|c|c||c|c|c}%\hline
		$J$ &$H$& $m$  & $e_{L^2}^u$   & $e_{a}^u$  \tabularnewline\hline
		6+0&$\sqrt{2}$/20	&3&7.367e-03 & 6.262e-02   \tabularnewline\hline 
		4+1&$\sqrt{2}$/20	&3&4.871e-03 & 4.186e-02 \tabularnewline\hline
		4+2&$\sqrt{2}$/20	&3&4.441e-03 & 4.018e-02  \tabularnewline%\hline
	\end{tabular}
	\caption{Numerical results with varying coarse-grid size $H$ for Test model $1$.} 
	\label{ta:nonlinear_online_h1}
\end{table}
\begin{table}
	\centering \begin{tabular}{c|c|c||c|c|c}%\hline
		$J$ &$H$& $m$  & $e_{L^2}^u$   & $e_{a}^u$  \tabularnewline\hline
		6+0&$\sqrt{2}$/20	&3&2.401e-03 & 4.469e-02   \tabularnewline\hline 
		4+1&$\sqrt{2}$/20	&3&4.295e-05 & 1.503e-03 \tabularnewline\hline
		4+2&$\sqrt{2}$/20	&3&1.656e-05 & 6.991e-04  \tabularnewline%\hline
	\end{tabular}
	\caption{Numerical results with varying coarse-grid size $H$ for Test model $2$.} 
	\label{ta:nonlinear_online_h2}
	%errors reduction using online basis
\end{table}

\subsection{Nonlinear poroelasticity case}\label{npc}

In this section, we present the numerical results of our method for solving
the nonlinear poroelasticity problems.
We set $\alpha=.9$, $M=10^6$.
The computational time $T := S\tau= 1$, and the time step size is chosen as $\tau := 1/20$.  The initial pressure is zero.

We will compare the CEM-GMsFEM solution with the fine-grid solution at the last time step $S$ (so that $S\tau = T$).  At the $(n+1)$th Picard iteration, to quantify the accuracy of our multiscale solutions, 
we use the following
relative weighted $\bfa{L}^2$ errors and energy errors:
\begin{align*}
e^u_{\bfa{L}^2}=\frac{||(\bfa{u}_{\tu{ms}}-\bfa{u}_{h})||_{\bfa{L}^2(\Omega)}}{||\bfa{u}_{h}||_{\bfa{L}^2(\Omega)}},\quad
e^u_{a}=\sqrt{\frac{a_n(\bfa{u}_{\tu{ms}}-\bfa{u}_{h},\bfa{u}_{\tu{ms}}-\bfa{u}_{h})}
	{a_n(\bfa{u}_{h},\bfa{u}_{h})}}\,,\\
e^p_{L^2}=\frac{||(p_{\tu{ms}}-p_{h})||_{L^2(\Omega)}}{||p_{h}||_{L^2(\Omega)}},\quad
e^p_{b}=\sqrt{\frac{b_n(p_{\tu{ms}}-p_{h},p_{\tu{ms}}-p_{h})}
	{b_n(p_{h},p_{h})}}\,,
\end{align*}
where the reference solution $(\boldsymbol{u}_{h}, p_{h})$ is computed via (\ref{ls})-(\ref{lp}) on the fine grid, the multiscale solution $(\bfa{u}_{\tu{ms}}, p_{\tu{ms}})$ is defined in (\ref{lsm})-(\ref{lpm}), and the bilinear forms $a_n$ and $b_n$ are defined in (\ref{biform1}) and (\ref{biform2}), respectively.

\begin{figure}
	\centering
	\subfigure{
    \includegraphics{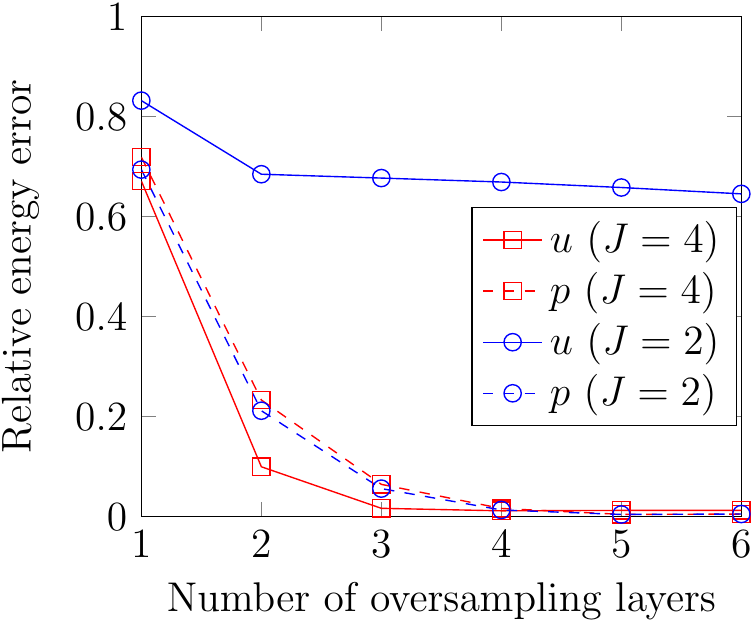}}
    \subfigure{\includegraphics{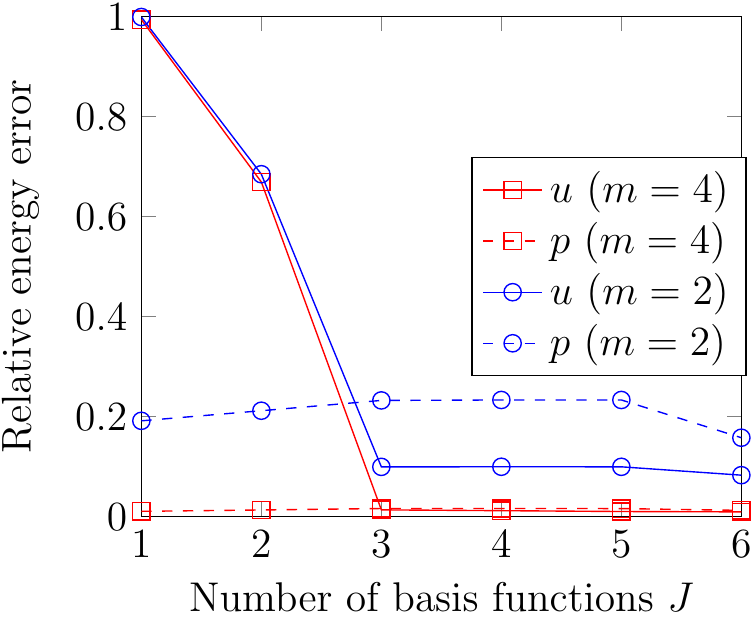}}
	\caption{Relative energy error (Test model 1) for $H=\sqrt{2}/20$ and fixed $J$ (left), fixed $m$ (right).}
	\label{fig:poro_ovnb1}
\end{figure}

\begin{figure}
	\centering
	\subfigure{
    \includegraphics{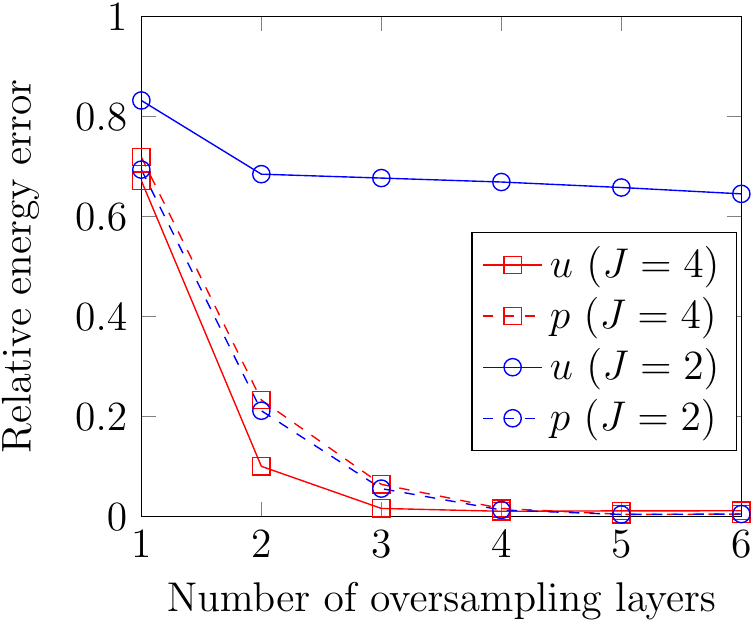}}
    \subfigure{\includegraphics{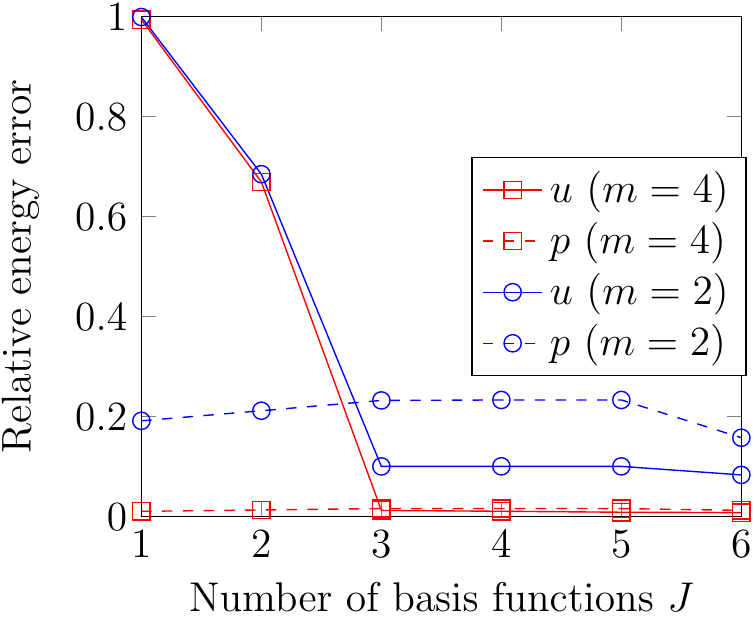}}
	\caption{Relative energy error (Test model 2) for $H=\sqrt{2}/20$ and fixed $J$ (left), fixed $m$ (right).}
	\label{fig:poro_ovnb2}
\end{figure}

We also first study the behavior of CEM-GMsFEM solution as $H$ becomes smaller.
The results are presented in Table \ref{ta:nonlinearporo_h1} and Table
\ref{ta:nonlinearporo_h2}. As expected, the accuracy of the CEM-GMsFEM solution
here improves for both the pressure and displacement as $H$ converges.
Figure \ref{fig:poro_ovnb1} and Figure \ref{fig:poro_ovnb2} display
the influence of the number of basis functions and oversampling layers.
Adding basis and oversampling layers will 
improve the accuracy of the displacement. We also find that more oversampling
layers yield more accurate pressure solution. However, the
accuracy of the pressure field is almost independent of 
the number of basis functions.  

\begin{table}
	\centering \begin{tabular}{c|c|c||c|c|c|c}%\hline
		$J$ &$H$& $m$  & $e_{L^2}^u$   & $e_{a}^u$& $e_{L^2}^p$   & $e_{b}^p$ \tabularnewline\hline
		4&$\sqrt{2}$/10	&3& 4.22e-03 &  3.05e-02 & 8.18e-04 & 1.82e-02     \tabularnewline\hline 
		4&$\sqrt{2}$/20	&4&2.09e-03 &  1.17e-02 & 5.08e-04 & 1.64e-02   \tabularnewline\hline
		4&$\sqrt{2}$/40	&5&2.06e-03 &  6.09e-03 & 2.82e-04 & 1.30e-02
		\tabularnewline%\hline
	\end{tabular}
	\caption{Numerical results with varying coarse-grid size $H$ for Test model $1$.} 
	\label{ta:nonlinearporo_h1}
\end{table}

\begin{table}
	\centering \begin{tabular}{c|c|c||c|c|c|c}%\hline
		$J$ &$H$& $m$  & $e_{L^2}^u$   & $e_{a}^u$& $e_{L^2}^p$   & $e_{b}^p$ \tabularnewline\hline
		4&$\sqrt{2}$/10	&3& 3.98e-03 &  3.05e-02 & 8.17e-04 & 1.82e-02   \tabularnewline\hline 
		4&$\sqrt{2}$/20	&4&8.24e-04 &  1.08e-02 & 5.07e-04 & 1.64e-02  \tabularnewline\hline
		4&$\sqrt{2}$/40	&5&2.59e-04 &  4.32e-03 & 2.79e-04 & 1.30e-02
		\tabularnewline
	\end{tabular}
	\caption{Numerical results with varying coarse-grid size $H$ for Test model $2$.} 
	\label{ta:nonlinearporo_h2}
	%errors reduction
\end{table}%%% %log(H)/log(1/10)*3

\vspace{50pt}
\section{Conclusions}\label{concs}
In this paper, we have proposed a framework of constraint energy minimizing generalized multiscale finite element method (CEM-GMsFEM) for solving problems of heterogeneous nonlinear poroelasticity (mainly) and elasticity.  In the case of nonlinear poroelasticity, the nonlinear stress equation involves strain-limiting elasticity and the nonlinear pressure equation is a Darcy-type parabolic equation.  Therefore, the key idea here is temporally discretizing the system by the implicit backward Euler scheme, then spatially linearizing it by Picard iteration (with desired convergence criterion) at each time step, until the terminal time.  In each Picard iteration, the CEM-GMsFEM is applied, to construct multiscale basis functions for both displacement and pressure systematically, with locally minimal energy, via using the techniques of oversampling, which leads to improved accuracy in the simulations.  Convergence analysis for each Picard iteration and numerical results has been shown to demonstrate the performance of the proposed method.  For the case of static nonlinear elasticity in the strain-limiting setting, the same strategy of Picard iteration combining with the CEM-GMsFEM is employed.  In addition to constructing the offline multiscale basis as in the poroelasticity case, we adaptively generate the residual based online basis, via solving a local problem in an oversampled domain with 
the residual as source.  Numerical tests prove the accuracy of our proposed method.  A proof of global convergence of the Picard iteration procedure is supplemented in Appendix \ref{conva}.

%%%%%

\vspace{20pt}

\noindent \textbf{Acknowledgements.}  

Eric Chung's work is partially supported by Hong Kong RGC General Research Fund
(Projects 14304217 and 14302018) and CUHK Direct Grant for Research 2018-19.  Tina Mai's work is funded by Vietnam National Foundation for Science and Technology Development (NAFOSTED) under grant number 101.99-2019.326.

%%%%%%%%%%%%%
\begin{appendix} 
\section{Comments on global convergence using Picard iteration algorithm, for nonlinear elasticity}\label{conva}

We will prove the global convergence of the Picard iteration procedure for our problem (\ref{4cem}) - (\ref{4cemn1}) by using fixed-point theorem, which \textbf{mainly} requires finding a suitable subset $\bfa{U}$ of the considering Banach space $\bfa{H}^1_0(\Omega)$, in which we are looking for solution $\bfa{u}$.  
%the main tasks to apply Fixed Point Theorems
%and related results is to identify a suitable subset of an appropriate function space in which, \cite{linear 1}
%we are looking for a solution

%This task might be worth a separate paper (if we can do, and it might be hard as well as time-consuming) because of our problem's nonlinearity?

In this paper, we will only introduce the key idea of the proof, where we can assume a Banach subspace $\bfa{U}$ of $\bfa{H}^1_0(\Omega)$ for $\bfa{u}$.  We want to find an operator $F: \bfa{u}^n \mapsto \bfa{u}^{n+1}$ such that $\bfa{u}^{n+1} = F(\bfa{u}^n), \|\bfa{F}'\|_*  <1\,, \bfa{F}' = DF(\bfa{u}^n)[\bfa{w}]
\in \mathcal{F}' \,,$ for any $\bfa{w} \in \bfa{H}^1_0(\Omega)$, for some suitable Banach space $\mathcal{F}' \subset \bfa{H}^1_0(\bfa{H}^1_0(\Omega))$ and some corresponding norm $\| \cdot \|_*$.  We have $\| \bfa{F}' \|_* \geq \| \bfa{F}' \|_1\,.$  Note that a mapping may be a contraction for some norm, but may not be a contraction for different norm.  Thus, identifying the right norm is important.

We denote $K(\bfa{u}^n):= \kappa(\bfa{x},|\bfa{Du}^n|)= \dfrac{1}{1-\beta(\bfa{x}) | \bfa{Du}^n|} $.  Given $\bfa{u}^n$, the next approximation $\bfa{u}^{n+1}$ is the solution of the system
\begin{align}
 -\ddiv (K( \bfa{u}^n) 
 \bfa{D}( \bfa{u}^{n+1})) &= \bfa{f} \quad \text{in } \Omega\,, \label{4bcem} \\
 \bfa{u}^{n+1} &= \bfa{0} \quad \text{on } \partial \Omega\, \label{4ccem}.
\end{align}

Now, for any $\bfa{w} \in \bfa{H}^1_0(\Omega)$ and very small $\e >0$, 
\begin{align*}
F(\bfa{u}^n + \e \bfa{w}) &= F(\bfa{u}^n) + DF(\bfa{u}^n) [\e \bfa{w}] + o(|\e \bfa{w}|)\\
&= \bfa{u}^{n+1} + DF(\bfa{u}^n) [\e \bfa{w}] + o(|\e \bfa{w}|)\,.
\end{align*}
%Taylor, and a function
%depending upon H that %tends to zero faster than |H|
Given $\bfa{u}^n + \e \bfa{w}$, let $\bfa{z}=\e DF(\bfa{u}^n)[\bfa{w}]$, where $DF(\bfa{u}^n)$ is a second order tensor.  Then, the next solution $F(\bfa{u}^n + \e \bfa{w})$ is the solution $\bfa{u}^{n+1}$ of the system 
\begin{align}
 -\ddiv (K( \bfa{u}^n + \e \bfa{w}) 
 \bfa{D}(\bfa{u}^{n+1} + \bfa{z})) &= \bfa{f} \quad \text{in } \Omega\,, \label{5bcem} \\
 \bfa{u}^{n+1} &= \bfa{0} \quad \text{on } \partial \Omega\, \label{5ccem}.
\end{align}

%%%
Subtracting (\ref{4bcem}) from (\ref{5bcem}), and dividing the result by $\e$, we obtain
\begin{align}\label{subs}
 -\ddiv \left(\frac{K( \bfa{u}^n + \e \bfa{w}) - K(\bfa{u}^n)}{\e}\bfa{Du}^{n+1}\right) 
 = \ddiv \left( K(\bfa{u}^n + \e \bfa{w}) \frac{\bfa{Dz}}{\e} \right)\,. 
\end{align}
Now, multiplying both sides of (\ref{subs}) by $\bfa{F}' = DF(\bfa{u}^n)[\bfa{w}]$, then letting $\e \to 0$, we get
\begin{align}\label{ibp1}
 \int_{\Omega} - \ddiv((DK(\bfa{u}^n)[\bfa{w}]) \bfa{Du}^{n+1}) \cdot \bfa{F}' \, \dx
 & = \int_{\Omega} \ddiv ( K(\bfa{u}^n) \bfa{DF}') \cdot \bfa{F}' \, \dx \,.
\end{align}
Integrating by parts both sides of (\ref{ibp1}), we obtain
\begin{align}\label{ibp2}
\int_{\Omega} (DK(\bfa{u}^n)[\bfa{w}]) (\bfa{Du}^{n+1}) \cdot (\bfa{DF}') \, \dx = - \int_{\Omega} K(\bfa{u}^n) |\bfa{DF}'|^2 \, \dx\,. 
 \end{align}

We assume that $K' = DK(\bfa{u}^n)[\bfa{w}] \in \mathcal{K}'$, for some suitable Banach space $\mathcal{K}'$, with some corresponding norm $\| \cdot \|$.  We have, $\| K' \| \leq \| K' \|_{\infty}$.  Also, from (\ref{ebound}), we note that
\begin{align}\label{ine7}
K(\bfa{u}^n) = \dfrac{1}{1-\beta(\bfa{x}) | \bfa{Du}^n|} >1\,.
\end{align}
Taking absolute values both sides of (\ref{ibp2}), then using inequality (\ref{ine7}) for the right hand side, and applying the Cauchy-Schwarz inequality to the left hand side of the result, we get
\begin{align}\label{ine7b}
\|DK(\bfa{u}^n)[\bfa{w}]\|_{\infty} \| \bfa{Du}^{n+1} \|_{\mathbb{L}^2(\Omega)} \|\bfa{DF}' \|_{\mathbb{L}^2(\Omega)} \geq \| \bfa{DF}' \|_{\mathbb{L}^2(\Omega)}^2 \,.
\end{align}
That is, 
\begin{align}\label{ine7c}
\|DK(\bfa{u}^n)[\bfa{w}]\|_{\infty} \| \bfa{Du}^{n+1} \|_{\mathbb{L}^2(\Omega)}  \geq \| \bfa{DF}' \|_{\mathbb{L}^2(\Omega)} \,.
\end{align}
Since $\| \bfa{f} \|_{\bfa{L}^2(\Omega)} \geq \| \bfa{Du}^{n+1} \|_{\mathbb{L}^2(\Omega)}$ (for the left hand side of (\ref{ine7c})),
%Poincare inequality for u, as take v = u in original equation
and $ \| \bfa{DF}' \|_{\mathbb{L}^2(\Omega)} \geq \| \bfa{F}' \|_1$ (for the right hand side of (\ref{ine7c})), it follows from (\ref{ine7c}) that
\[ \| K' \|_{\infty} \| \bfa{f} \|_{\bfa{L}^2} > \| \bfa{F}' \|_1\,.\]
We can choose $\bfa{f}$ at the beginning such that $\| \bfa{f} \|_{\bfa{L}^2}$ can dominate $\| K' \|_{\infty}$ in the way that $1 > \| K' \|_{\infty} \| \bfa{f} \|_{\bfa{L}^2} \; (> \| \bfa{F}' \|_1)$, and we are done.

To find the expression of 
$K' = DK(\bfa{u}^n)[\bfa{w}] \in \mathcal{K}'$, we compute (by definition of the Fr\'{e}chet derivative) as follows.  In preparation, let $\phi(\bfa{A}) = |\bfa{A}|$.  Then, by Taylor expansion, we get
\begin{align}\label{dernorm}
\begin{split}
 | \bfa{Du}^n + \bfa{D}(\e \bfa{w}) | &=\phi(\bfa{Du}^n + \bfa{D}(\e \bfa{w})) \\
 &= |\bfa{Du}^n + \e \bfa{D}( \bfa{w})|\\
 &= |\bfa{Du}^n| + D_{(\bfa{Du}^n)}(|\bfa{Du}^n|) [\bfa{D}(\e \bfa{w})] + o(\e(\bfa{Dw}))\\
 & = |\bfa{Du}^n| + \frac{(\bfa{Du}^n) \cdot (\bfa{D}(\e \bfa{w}))}{|\bfa{Du}^n|} + o(\e\bfa{Dw})\,,
 \end{split}
\end{align}
where, in the last equality, we use the result
\[D_{\bfa{A}}(|\bfa{A}|) [\bfa{W}] = D_{\bfa{A}}(\bfa{A} \cdot \bfa{A})^{1/2} [\bfa{W}] = \dfrac{1}{2} (\bfa{A} \cdot \bfa{A})^{-1/2} D_{\bfa{A}}(\bfa{A} \cdot \bfa{A})[\bfa{W}] = \dfrac{\bfa{W} \cdot \bfa{A}}{|\bfa{A}|}\,.\]
By the Riesz Representation Theorem, for $DK(\bfa{u}^n)[\bfa{w}]$, there exists a unique element, namely Gradient of $K(\cdot)$ at $\bfa{u}^n$, denoted by $\nabla K (\bfa{u}^n) \in \bfa{H}^1_0(\Omega)$ such that
\begin{align}\label{rrep}
DK(\bfa{u}^n)[\bfa{w}] = \nabla K (\bfa{u}^n) \cdot \bfa{w}\,, \forall \bfa{w} \in \bfa{H}^1_0(\Omega)\,.
\end{align}
%%%
To find $\nabla K (\bfa{u}^n)$, we compute as follows:
\begin{align*}
 K' &= DK(\bfa{u}^n)[\bfa{w}]\\
 &= \lim_{\e \to 0} \frac{K(\bfa{u}^n + \e w) - K(\bfa{u}^n)}{\e}\\
 &= \lim_{\e \to 0} \frac{\dfrac{1}{1 - \beta(\bfa{x}) | \bfa{D}(\bfa{u}^n + \e \bfa{w})|} - \dfrac{1}{1- \beta(\bfa{x})|\bfa{Du}^n|}}{\e}\\
 &= \lim_{\e \to 0} \frac{\beta(\bfa{x})(|\bfa{D}(\bfa{u}^n+ \e \bfa{w})| - |\bfa{Du}^n|)}{\e(1 - \beta(\bfa{x}) | \bfa{D}(\bfa{u}^n + \e \bfa{w})|)(1- \beta(\bfa{x})|\bfa{Du}^n|)}\\
 &= \lim_{\e \to 0}\frac{\beta(\bfa{x}) \e \dfrac{(\bfa{Du}^n) \cdot (\bfa{Dw})}{|\bfa{Du}^n|} + o (\e\bfa{Dw})}{\e(1 - \beta(\bfa{x}) | \bfa{D}(\bfa{u}^n + \e \bfa{w})|)(1- \beta(\bfa{x})|\bfa{Du}^n|)}\\
 &= \frac{\beta(\bfa{x})(\bfa{Du}^n) \cdot (\bfa{Dw})}{|\bfa{Du}^n|(1- \beta(\bfa{x})|\bfa{Du}^n|)^2}\,,
\end{align*}
where the last expression follows from (\ref{dernorm}).  Integrating the last expression by parts, we get
\[\int_{\Omega}\frac{\beta(\bfa{x})(\bfa{Du}^n) \cdot (\bfa{Dw})}{|\bfa{Du}^n|(1- \beta(\bfa{x})|\bfa{Du}^n|)^2} \, \dx =\int_{\Omega}
\beta(\bfa{x}) \ddiv \left( \dfrac{\bfa{Du}^n}{|\bfa{Du}^n|(1- \beta(\bfa{x})|\bfa{Du}^n|)^2} \right) \cdot \bfa{w} \, \dx\,,\]
in which
\begin{align}\label{gradient}
\nabla K(\bfa{u}^n) = \beta(\bfa{x}) \ddiv \left( \dfrac{\bfa{Du}^n}{|\bfa{Du}^n|(1- \beta(\bfa{x})|\bfa{Du}^n|)^2} \right)\,,
 \end{align}
 as expected.
From (\ref{rrep}), we can assume that $\nabla K(\bfa{u}^n) \in \mathcal{H}$, for some Banach subspace $\mathcal{H}$ of $\bfa{H}^1_0(\Omega)$, with norm $\| \cdot \|_{\mathcal{H}}$.  It holds that $\| \nabla K(\bfa{u}^n) \|_1 \leq \| \nabla K(\bfa{u}^n) \|_{\mathcal{H}}$.  Currently, we have not known the exact Banach subspace $\mathcal{H}$ even we know $\bfa{Du}^n \in \mathbb{L}^{\infty}(\Omega)$.  The reason lies in the denominator of (\ref{gradient}):  From (\ref{calZ}), there is an upper bound of $|\bfa{Du}^n|$; but we do not know whether it has a maximum.  (If there is a $\bfa{u}^n$ such that the maximum of $|\bfa{Du}^n|$ is attained, then we do not know whether such a $\bfa{u}^n$ satisfies the boundary value problem (\ref{form3}) - (\ref{D}).)

\end{appendix}

%%%%%%%%%%%
\bibliographystyle{plain} 

\bibliography{1m_elliptic1}

\begin{thebibliography}{10}

\bibitem{unitybabu}
I.~Babuska and J.~M. Melenk.
\newblock The partition of unity method.
\newblock {\em International Journal of Numerical Methods in Engineering},
  40:727--758, 1996.

\bibitem{Beck2017}
Lisa Beck, Miroslav Bul{\'i}{\v{c}}ek, Josef M{\'a}lek, and Endre S{\"u}li.
\newblock On the existence of integrable solutions to nonlinear elliptic
  systems and variational problems with linear growth.
\newblock {\em Archive for Rational Mechanics and Analysis}, 225(2):717--769,
  Aug 2017.

\bibitem{biot}
M.~A. {Biot}.
\newblock {General theory of three-dimensional consolidation}.
\newblock {\em Journal of Applied Physics}, 12:155--164, February 1941.

\bibitem{nlss}
Michele Botti, Daniele~A. Di~Pietro, and Pierre Sochala.
\newblock A nonconforming high-order method for nonlinear poroelasticity.
\newblock In Cl{\'e}ment Canc{\`e}s and Pascal Omnes, editors, {\em Finite
  Volumes for Complex Applications VIII - Hyperbolic, Elliptic and Parabolic
  Problems}, pages 537--545, Cham, 2017. Springer International Publishing.

\bibitem{dlp}
Donald~L. Brown and Maria Vasilyeva.
\newblock {A Generalized Multiscale Finite Element Method for poroelasticity
  problems I: Linear problems}.
\newblock {\em Journal of Computational and Applied Mathematics}, 294:372 --
  388, 2016.

\bibitem{gnone}
Donald~L. Brown and Maria Vasilyeva.
\newblock A generalized multiscale finite element method for poroelasticity
  problems {II}: Nonlinear coupling.
\newblock {\em Journal of Computational and Applied Mathematics}, 297:132 --
  146, 2016.

\bibitem{B-M-S}
M.~Bul{\'i}{\u c}ek, J.~M{\'a}lek, and E.~S{\"u}li.
\newblock {Analysis and approximation of a strain-limiting nonlinear elastic
  model}.
\newblock {\em Mathematics and Mechanics of Solids}, 20(I):{92--118}, 2015.
\newblock DOI: 10.1177/1081286514543601.

\bibitem{BMRS14}
Miroslav Bul\'{i}\v{c}ek, Josef M\'{a}lek, K.~R. Rajagopal, and Endre S\"{u}li.
\newblock {On elastic solids with limiting small strain: modelling and
  analysis}.
\newblock {\em EMS Surveys in Mathematical Sciences}, 1(2):283--332, 2014.

\bibitem{chung2016adaptive}
Eric Chung, Yalchin Efendiev, and Thomas~Y Hou.
\newblock Adaptive multiscale model reduction with generalized multiscale
  finite element methods.
\newblock {\em Journal of Computational Physics}, 320:69--95, 2016.

\bibitem{chungres}
Eric~T. Chung, Yalchin Efendiev, and Wing~Tat Leung.
\newblock Residual-driven online generalized multiscale finite element methods.
\newblock {\em Journal of Computational Physics}, 302:176 -- 190, 2015.

\bibitem{cem1}
Eric~T. Chung, Yalchin Efendiev, and Wing~Tat Leung.
\newblock Constraint energy minimizing generalized multiscale finite element
  method.
\newblock {\em Computer Methods in Applied Mechanics and Engineering}, 339:298
  -- 319, 2018.

\bibitem{cem2f}
Eric~T. Chung, Yalchin Efendiev, and Wing~Tat Leung.
\newblock Fast online generalized multiscale finite element method using
  constraint energy minimization.
\newblock {\em Journal of Computational Physics}, 355:450 -- 463, 2018.

\bibitem{chungres1}
Eric~T. Chung, Yalchin Efendiev, and Guanglian Li.
\newblock {An adaptive GMsFEM for high-contrast flow problems}.
\newblock {\em Journal of Computational Physics}, 273:54 -- 76, 2014.

\bibitem{C-G-K}
P.~G. Ciarlet, G.~Geymonat, and F.~Krasucki.
\newblock {A new duality approach to elasticity}.
\newblock {\em Mathematical Models and Methods in Applied Sciences}, 22(1):{21
  pages}, 2012.
\newblock DOI: 10.1142/S0218202512005861.

\bibitem{G1}
Yalchin Efendiev, Juan Galvis, and Thomas~Y. Hou.
\newblock {Generalized multiscale finite element methods (GMsFEM)}.
\newblock {\em J. Comput. Phys.}, 251:116--135, October 2013.

\bibitem{G2}
Yalchin Efendiev, Juan Galvis, Guanglian Li, and Michael Presho.
\newblock {Generalized multiscale finite element methods. Nonlinear elliptic
  equations}.
\newblock {\em Communications in Computational Physics}, 15(3):733--755, 003
  2014.

\bibitem{lporo}
Shubin Fu, Robert Altmann, Eric~T. Chung, Roland Maier, Daniel Peterseim, and
  Sai-Mang Pun.
\newblock Computational multiscale methods for linear poroelasticity with high
  contrast.
\newblock {\em Journal of Computational Physics}, 395:286 -- 297, 2019.

\bibitem{gne}
Shubin Fu, Eric Chung, and Tina Mai.
\newblock Generalized multiscale finite element method for a strain-limiting
  nonlinear elasticity model.
\newblock {\em Journal of Computational and Applied Mathematics}, 359:153 --
  165, 2019.

\bibitem{cemgle}
Shubin {Fu} and Eric~T. {Chung}.
\newblock {Constraint energy minimizing generalized multiscale finite element
  method for high-contrast linear elasticity problem}.
\newblock {\em Accepted by Communications in Computational Physics. arXiv
  e-prints}, page arXiv:1809.03726, Sep 2018.

\bibitem{B-Mai-Walton}
Tina Mai and Jay~R. Walton.
\newblock {On monotonicity for strain-limiting theories of elasticity}.
\newblock {\em Journal of Elasticity}, 120(I):{39--65}, 2015.
\newblock DOI: 10.1007/s10659-014-9503-4.

\bibitem{A-Mai-Walton}
Tina Mai and Jay~R. Walton.
\newblock {On strong ellipticity for implicit and strain-limiting theories of
  elasticity}.
\newblock {\em Mathematics and Mechanics of Solids}, 20(II):{121--139}, 2015.
\newblock DOI: 10.1177/1081286514544254.

\bibitem{korn}
Patrizio Neff, Dirk Pauly, and Karl-Josef Witsch.
\newblock {Poincar\'{e} meets Korn via Maxwell: Extending Korn's first
  inequality to incompatible tensor fields}.
\newblock {\em Journal of Differential Equations}, 258(4):1267 -- 1302, 2015.

\bibitem{Raji03}
K.~R. Rajagopal.
\newblock On implicit constitutive theories.
\newblock {\em Applications of Mathematics}, 48(4):279--319, 2003.

\bibitem{KRR-ZAMP2007}
K.~R. Rajagopal.
\newblock {The elasticity of elasticity}.
\newblock {\em Z. Angew. Math. Phys.}, {58}({2}):{309--317}, {2007}.

\bibitem{KRR-MMS2011b}
K.~R. Rajagopal.
\newblock {Conspectus of concepts of elasticity}.
\newblock {\em Mathematics and Mechanics of Solids}, {16}({5, SI}):{536--562},
  {2011}.

\bibitem{KRR-MMS2011a}
K.~R. Rajagopal.
\newblock {Non-linear elastic bodies exhibiting limiting small strain}.
\newblock {\em Mathematics and Mechanics of Solids}, {16}({1}):{122--139},
  {2011}.

\bibitem{KRR-ARS-PRSA2007}
K.~R. Rajagopal and A.~R. Srinivasa.
\newblock {On the response of non-dissipative solids}.
\newblock {\em Proceedings of the Royal Society of London, Mathematical,
  Physical and Engineering Sciences}, {463}({2078}):{357--367}, {2007}.

\bibitem{Rajagopal493}
K.R Rajagopal and A.R Srinivasa.
\newblock On a class of non-dissipative materials that are not hyperelastic.
\newblock {\em Proceedings of the Royal Society of London A: Mathematical,
  Physical and Engineering Sciences}, 465(2102):493--500, 2009.

\bibitem{eulp}
R.E. Showalter.
\newblock Diffusion in poro-elastic media.
\newblock {\em Journal of Mathematical Analysis and Applications}, 251(1):310
  -- 340, 2000.

\end{thebibliography}
\end{document}